%% file: SPAJamming15.tex
\documentclass[3p]{elsarticle}

\usepackage{graphicx}

\usepackage{enumerate}
\usepackage{latexsym,amssymb,amsthm,amsxtra}
\usepackage{amsmath,mathrsfs}
\usepackage{tikz}
\usepackage{pgf}
\usepackage{dcolumn}
\usepackage{stmaryrd}
\usepackage{pstricks}
\usepackage{pst-node}
\usepackage{multido}

\def\ep{\hfill $\Box$}

\def\bp{\noindent{\bf Proof.}\ }
\def\epsilon{\varepsilon}

\usepackage{anysize}
\marginsize{3cm}{2.5cm}{2.5cm}{2.5cm} 
\newtheorem{theorem}{Theorem}[section]
\newtheorem{lemma}[theorem]{Lemma}
\newtheorem{prop}[theorem]{Proposition}

\newtheorem{coro}[theorem]{Corollary}
\newtheorem{remark}[theorem]{Remark}
\newtheorem{definition}[theorem]{Definition}

\newtheorem{example}[theorem]{Example}

\newcommand{\be}{ \begin{equation}}
\newcommand{\ee}{\end{equation}}

\def\E{{\mathbb E}}

\def\P{{\mathbb P}}

\def\R{{\mathbb R}}
\def\Z{{\mathbb Z}}
\def\V{{\cal V}}
\def\N{{\mathbb N}}

\def\A{{\mathscr{A}}}
\def\B{{\mathscr{B}}}
\def\maC{{\mathscr{C}}}
\def\D{{\mathcal{D}}}
\def\F{{\mathcal{F}}}
\def\T{{\mathcal{T}}}

\def\barmun{\bar {\mu}^n}
\def\barmunt{\bar{\mu}^{n,\tau^n_{\alpha}}}
\def\barmu{\bar\mu}
\def\barmus{\bar\mu^*}
\def\sA{\mbox{\textsc{a}}}
\def\sB{\mbox{\textsc{b}}}
\def\sU{\mbox{\textsc{u}}}
\def\V{{\mathcal V}}

\def\maU{{\mathcal U}}
\def\maA{{\mathcal{A}}}
\def\maV{{\mathcal{V}}}
\def\maB{{\mathcal{B}}}

\def\maH{{\mathcal{H}}}

\def\C{{\mathcal{C}}}

\def\T2a{{\tau_{2\alpha}}}
\def\t2a{{t_{2\alpha}}}
\newcommand\cro[1]{\left\langle #1 \right\rangle}
\newcommand\tend{\underset{n \to \infty}{\longrightarrow}}
\newcommand\tendP{\overset{(\mathcal P)}{\underset{n \to \infty}{\longrightarrow}}}
\newcommand\proc[1]{\left(#1_t\right)_{t \ge 0}}

\newcommand\suite[1]{\left\{#1\right\}_{n\in\N^*}}

\def\M{{\cal{M}\rm}}

\def\({{\Bigl(}}
\def\){{\Bigr)}}
\newcommand\pr[1]{{\mathbf P}\left[#1\right]}
\newcommand\esp[1]{{\mathbf E}\left[#1\right]}

\def\ind{{\mathchoice {\rm 1\mskip-4mu l} {\rm 1\mskip-4mu l}
{\rm 1\mskip-4.5mu l} {\rm 1\mskip-5mu l}}}
\def\square{\ifmmode\sqr\else{$\sqr$}\fi}
\def\sqr{\vcenter{
         \hrule height.1mm
         \hbox{\vrule width.1mm height2.2mm\kern2.18mm\vrule width.1mm}
         \hrule height.1mm}}                  

\title{The Jamming Constant of Uniform Random Graphs}

\begin{document}

\begin{frontmatter}

\title{The Jamming Constant of Uniform Random Graphs}


\author[mymainaddress]{Paola Bermolen}
\ead{paola@fing.edu.uy}

\author[mysecondaryaddress]{Matthieu Jonckheere}
\ead{matthieu.jonckheere@gmail.com}

\author[mythirdaddress]{Pascal Moyal\corref{mycorrespondingauthor}}
\cortext[mycorrespondingauthor]{Corresponding author}
\ead{pascal.moyal@utc.fr}

\address[mymainaddress]{Universidad de La Republica, Av. De Julio 1968, Montevideo, Uruguay}
\address[mysecondaryaddress]{CONICET, Mathematics Department, Facultad de Ciencas Exactas y Naturales, Universidad de Buenos Aires, 1482 Pabellon I, Ciudad Universitaria Buenos Aires, Argentina}
\address[mythirdaddress]{Laboratoire de Mathematiques Appliquees LMAC, Universite de Technologie de Compiegne, Rue du Dr Schweitzer
CS 60319, 60203 Compiegne, France}

\begin{abstract}
By constructing jointly a random graph and an associated exploration process,
we define the dynamics of a ``parking process'' on a class of uniform random graphs as a measure-valued Markov process, 
representing the empirical degree distribution of non-explored nodes.
We then establish a functional law of large numbers for this process as the number of vertices grows to infinity, allowing us to assess 
the jamming constant of the considered random graphs, i.e. the size of the maximal independent set discovered by the exploration algorithm. 
This technique, which can be applied to any uniform random graph with a given degree distribution, can be seen as a generalization in the space of measures, 
of the differential equation method introduced by Wormald.
\end{abstract}



\end{frontmatter}


\section{Introduction}

Consider a finite graph $G$ for which $\V$ is the set of nodes or sites.
The parking process in continuous time on $G$ may be described as follows.
At time $0$, all sites are vacant. They all have independent exponential clocks.
When the clock of a given vacant site rings and all of its neighbors are vacant, this site turns occupied. Otherwise, nothing happens.
Once occupied, a site remains so for ever.
The process goes on until all sites are either occupied or have at least one of their neighbors occupied.
The final state of the process is often referred to as the jamming limit of $G$, and the final proportion of 
occupied sites, its jamming constant. 

Our motivation to study the parking process on random graphs is twofold.
On the one hand, these dynamics are the simplest procedure to discover maximal independent sets and
have been extensively studied for some specific graphs. 
Explicit results have been obtained for regular graphs  \cite{wormaldDF}, exploiting 
 their very specific structure (see also \cite{gamarnik} for graphs with large girths). 
In the Erd\"os R\'enyi case, a similar differential method can be employed thanks to the great amount of independence and symmetry of the collection of edges, to get an explicit jamming constant (see Theorem 2.2 and the references in \cite{mcdiarmid}). 
Hence, to look at ``uniform" random graphs having a given asymptotic degree distribution, but much less structure and symmetry, 
is a natural continuation of this research avenue.

On the other hand, parking processes have received a great amount of attention in the case of spatial structures.
It has been considered on discrete structures like $\mathbb Z^d$ \cite{ritchie,ferrari-fernandez-garcia-02}
and on point processes \cite{penrose,ferrari-fernandez-garcia-02,baccelli-09,baccelli-12}. 
In physics and biological sciences, where it is usually referred to as random sequential absorption, it models phenomena of deposition of colloidal particles or proteins on surfaces (see \cite{evans1993}).
In communication sciences and in wireless networks in particular,
it allows to represent the number of connections for CSMA-like algorithms in a given time-slot,
for a given spatial configuration of terminals (see \cite{kleinrock} for a classical reference on the definition of the protocol).
The general idea of CSMA is to schedule transmissions in such a way that nodes that interfere each other would not transmit simultaneously, 
see for instance \cite{baccelli-09} for a stochastic geometry-based model in which CSMA is approximated by a Mat\'ern-like process.
Unfortunately, spatial models are in general very difficult to study theoretically and to the best of our knowledge,
there is no efficient way to compute the jamming constant in most cases.
Studying the jamming constant of uniformly chosen random graphs with a given asymptotic degree distribution 
(we make this notion more precise in the sequel) allows to make a first step in this direction, by studying a ``first order" model which grasps only
the bonds between points but no further correlations.
The techniques and analysis presented here are at the core of the performance evaluation analysis of wireless systems, as developed in a companion applied paper \cite{bermolen}.

In this paper we focus on the parking process for a class of random graphs having given deterministic 
asymptotic degree distributions, and derive a computable characterization of the jamming constant, as the number of vertices grows to infinity. 
To describe the evolution of the parking process in a Markovian fashion, without keeping track of a too large set of information, our strategy is the following: we start from the degree distribution of the graph, and then construct simultaneously the random graph and the associated parking process.
In both cases, the underlying (multi-)graph is constructed similarly to the configuration model,
 see \cite{wormaldCM,remco,janson2014} and the references therein. 
A similar approach was considered in \cite{DDMT12}, to construct a random social network together with a SIR process which propagates on it.

Following these ideas, we first define these dynamics for a graph having a fixed number of nodes $n$ and study the time-evolution of the empirical measure of the degrees of the vacant sites, which defines a measure-valued Markov process. 
Then, under the assumption that the initial empirical measure of degrees converges to a measure having mild moment assumptions, 
we take $n$ to infinity and prove a functional law of large numbers on the evolution of the empirical measures of degrees.
We show in particular that given our assumptions on the initial random degree distribution, 
the limit is unique, and defined as the solution  of a non-linear infinite-dimensional system of differential equations.
This can be seen as a generalization in the space of point measures, of the differential equation method introduced by Wormald \cite{wormaldDF}, 
providing simple and insightful characteristics of the random graphs under consideration, as solutions of finite-dimensional differential equations.   
To the best of our knowledge, this is the first such general result, 
embracing in particular, several particular cases investigated in the literature.


In the case of the Poisson distribution, we are able to explicitly compute the measure-valued flow of unexplored nodes, which
turns out to be an inhomogeneous Poisson measure. We then retrieve the jamming constant of the Erd\"os-R\'enyi (ER) graph. 
(Note however that our construction does not lead to a proper ER graph, only to a random multi-graph having the same degree distribution.) 
We also retrieve a constant calculated by R\'enyi \cite{finch, solomon} for a spatial model on $ \mathbb Z$, showing that both models share the same jamming limit.

The proof of the functional law of large numbers is based on successive approximations for the generator of the infinite-dimensional Markov process, relying in particular on quantifying the probability to obtain self-loops and multi-edges. Note that this difficulty is inherent to
 configuration model constructions, which have the disadvantage of constructing a multi-graph rather than a simple one, though elegant arguments have shown that with a probability independent of the size of the graph,  a simple graph is obtained \cite{janson2014} (see also the monograph \cite{remco}).
The uniqueness of the deterministic limiting measure-valued flow is not immediate and has to be proved using an adequate norm on the spaces of solutions.

The rest of the paper is organized as follows.
 In Section \ref{sec:construct}, we describe the simultaneous construction of the parking process
and the random graph. In Section \ref{sec:generator} we calculate the generator of the induced measure-valued Markov process, and the corresponding 
 semi-martingale decomposition is introduced. 
In Section \ref{sec:hyd}, we state our main result and its consequences. In particular we show how the latter leads to closed-forms, or at least to computable characterizations of the jamming constant in various cases. 
Section \ref{sec:proof} is devoted to the proofs of our main results.


\newpage.

\paragraph{Notation}


Let us introduce the main notation used throughout the paper.
\begin{itemize}
\item We denote by $\R$ the set of real numbers, and $\R+$ (respectively, $\R^*$) the subset of non-negative (resp., non-null) real numbers.
Let also $\N$ be the set of non-negative integers and $\N^*$, the subset of positive integers. For any $x,y \in \R$, let
$x\wedge y=\min\{x,y\}$, $x \vee y=\max\{x,y\}$ and $x^+=x\vee 0$. Let also for any $a,b \in \N$, $\llbracket a,b \rrbracket=\{a,a+1,...,b\}.$ 
We denote by $o(.)$, a function: $\N \to \R$ such that $n^\beta o(n) \tend 0$ for any real number $\beta<1$. 
\item Let $\mathcal B_b$ be the set of Borel bounded functions: $\R \to \R$. For all $\phi \in\mathcal B_b$, denote
$$\parallel \phi \parallel=\sup_{x\in\R} \mid \phi(x) \mid.$$
Denote for any Borel set $A$, $\ind_A$ the indicator function of $A$. Denote by $\mathbf{1}$, the real function constantly equal to 1 and for any  $k \in \N$, $\chi^k$ the function $x \mapsto x^k$.
For all $\phi\in \mathcal B_b$, we also denote by $\Delta\phi$ the {\em discrete gradient} of $\phi$, \emph{i.e.} 
$$\Delta\phi (i)=\phi(i)-\phi(i-1),\,\forall\, i\in\N^*.$$
\item Let $\M_F(\N)$ be the set of finite measures on $\N$. We write $\mu(i):=\mu(\{i\})$ for any $\mu \in \M_F(\N)$ and any $i\in\N$.
The null measure is denoted $\mathbf 0$.
For all $\mu \in \M_F(\N)$ and all $\phi:\R \to \R$, $\cro{\mu,\phi} $ denotes
 the integral of $\phi$ with respect to $\mu$:
\begin{equation*}
\cro{\mu,\phi} = \int \phi(x) \mu(dx)=\sum\limits_{i\in\N} \phi(i)\mu(i).
\end{equation*}
In this way, for any such $\mu$ and any $A \subset \N$, $\cro{\mu,\ind_A}=\mu(A)$ is the measure of $A$,
$\cro{\mu,\mathbf 1}=\mu(\N)$ is the total mass of $\mu$, and for any $k \in \N^*$, $\cro{\mu, \chi^k}$ is the $k$-th moment of $\mu$.
For any $\mu\in\M_F(\N)$, denote $F_{\mu}$ the cumulative function associated to $\mu$ and $F^{-1}_\mu$, its generalized inverse.
For any counting measure $\mu \in \M_F(\N)$, we will be led to order and index the atoms of $\mu$ as follows:
\begin{itemize}
\item we denote for any $\ell \in \{1,...,\mu(0)\}$, $v_{\ell}(\mu)$ the $\ell$th atom of degree 0 ranked in arbitrary order;\\
\item by induction, for any  $i\in\N$ and any $\ell \in \llbracket 1,\mu(i+1)\rrbracket$, $v_{\sum_{j=1}^i \mu(j)+\ell}\left(\mu\right)$, the
$\ell$th atom of degree $i+1$, in arbitrary order, in a way that
\end{itemize}
\begin{equation}
\label{eq:notationv}
\mu=\sum_{j=1}^{\cro{\mu,\mathbf 1}} \delta_{{v_j}\left(\mu\right)}.
\end{equation}

\item For any Polish space $E$, $\D\left([0,T], E\right)$ is the Polish space of rcll functions from $[0,T]$ to $E$, and
$\C\left([0,T], E\right)$ is the space of continuous functions from $[0,T]$ to $E$. Both $\D\left([0,T], E\right)$ and
$\C\left([0,T], E\right)$ are equipped with the Skorokhod $J_1$-topology.
\item Unless explicitly mentioned, throughout all the random variables (r.v.'s, for short) are defined on a common probability space
$\left(\Omega,\mathcal F,\mathbf P\right).$ On the latter, let us write ``$\Rightarrow$" for weak convergence of r.v.'s, and 
``$\overset{(\mathcal P)}{\rightarrow}$" for convergence in probability. Finally, let us denote
$\proc{\langle\!\langle M \rangle\!\rangle}$ the quadratic variation of the $\D\left([0,T],E\right)$-valued martingale $\proc{M}$.
\end{itemize}

\section{Construction of the graph and Markov representation}
\label{sec:construct}

In this section we present our construction of a random graph of prescribed degree distribution, simultaneously with the parking process on the latter graph. The basic objects of our construction are:
\begin{itemize}
\item[(i)] a probability measure $\nu$ on $\N$ having support $\llbracket 0,n-1 \rrbracket$ (where $n\ge 1$), which will be referred to as the 
{\em degree distribution};
\item[(ii)] a $n$-independent sample $\mathbf d:=\left(d(1),...,d(n)\right)$ of the distribution $\nu$, termed {\em degree vector}.\footnote{We assume independence of the degrees for simplicity, however it should be noted that the results
hereafter hold in larger generality: for assumption (\ref{eq:convinit}) to hold true, we only need the convergence of the Ces\`aro means of the vector 
$\mathbf d$ as $n$ goes large.}  
The {\em empirical degree distribution} is the following random point measure,  
$$\mu_0=\sum_{i=1}^{n} \delta_{d(i)};$$   
\item[(iii)] the set of {\em nodes} $\maV$, whose elements are denoted $u_0(1),...,u_0(n)$ (the use of this notation will become clear in a few lines). 
We set a one-to-one relation between the nodes and the atoms 
of $\mu_0$ as follows: to the node $u_0(i)$ is associated the element $d\left(\gamma(i)\right)$, where $\gamma$ is a permutation of 
$\llbracket 1,n \rrbracket$ arranging $d$ in increasing order. We then say that $d(\gamma(i))$ is the {\em degree} of $u_0(i)$ 
and we write $d_{u_0(i)}:=d(\gamma(i))$. 
\end{itemize}

At time 0, all the nodes of $\maV$ are disconnected: the {\em associated graph} of our construction at time 0, denoted $G_0$, thus consists in the set 
of nodes $\maV$, without any edge. At this point, the nodes are all said {\em unexplored} (we say that they are ``$\sU$-nodes").  
We consider that each node has as many {\bf unmatched half-edges} as its degree - we then say that the node is the {\em ego} of its half-edges. 

Let us define for all $t \ge 0$, $\maU_t$, $\maA_t$ and $\maB_t$ the sets of {\em unexplored}, {\em active} and {\em blocked} nodes, respectively. At $t=0$, we thus fix $\maU_0=\maV$ (hence the notation above), and set $\maA_0=\maB_0=\emptyset$. Let us also define for all $t\ge 0$, 
\begin{equation*}
\maH_t= \Bigl\{\mbox{unmatched half-edges at time $t$}\Bigl\}.
\end{equation*}

 Let for any $t$ and any $j \in \maU_t$, $d_j(\maU_t)$ denote the number of unmatched half-edges of $j$ at time $t$. 
 Define also the following element of $\M_F(\N)$,  
 \begin{equation*}
 \mu_t=\sum_{j\in \maU_t} \delta_{d_j(\maU_t)},
 \end{equation*}
 termed {\em empirical degree distribution at }$t$. Notice that the cardinality $U_t$ of $\maU_t$, and the cardinality $H_t$ of $\maH_t$ can respectively be retrieved from $\mu_t$ by 
 \begin{gather*}
 U_t=\cro{\mu_t,\mathbf 1}\quad \mbox{ and } \quad H_t=\cro{\mu_t,\chi}.
 \end{gather*}
All these time-dependent quantities will be updated, by induction on the event times, as will be described hereafter. 

Fix $\lambda>0$, and let $\xi_0$ be a random variable of exponential distribution of rate $\lambda n$. 
As long as this exponential clock does not ring, the system remains unchanged: we set $\mu_t=\mu_0$ and likewise, 
$\maU_t,\maA_t,\maB_t$ and $\maH_t$ equal their initial value, for any $t\in [0,\xi_0)$. 
The dynamics of the system is then determined by induction, as follows:  
assume that a clock rings at time $t$. Then, several state changes occur instantaneously at $t$, but following a given sequence. 
To represent these changes, all characteristics of the system just before the event 
occurs are indexed by $t^-$ and then, to distinguish between the different steps of the state actualization, by $t^{-+}$, $t^{-++}$ 
and finally $t$.  

\begin{itemize}
\item[{\bf Step 1.}] \textbf{A node $a$ becomes active}: we draw uniformly an index $i$ in $\llbracket 1,  U_{t^-} \rrbracket$. 
The corresponding 
node $a:=u_{t^-}(i)$ of  $\maU_{t^-}$ becomes active and we set 
$$\maU_{t^{-+}}=\maU_{t^-} \setminus \{a\},\quad\maA_{t^{-+}}=\maA_{t^-} \cup \{a\}\quad\mbox{ and }\quad\maB_{t^{-+}}=\maB_{t^-}.$$
Let $K\left(\mu_{t^-}\right)=d_a$, the number of unmatched half-edges of $a$ at $t^-$. As $a$ is no longer a $\sU$-node, we update the measure 
$\mu_{t^-}$ as follows, 
\begin{equation}
\label{eq:genetape1}
\mu_{t^{-+}}=\mu_{t^-}-\delta_{K\left(\mu_{t^-}\right)}.
\end{equation}

\begin{figure}
\centering
\includegraphics[scale=0.9]{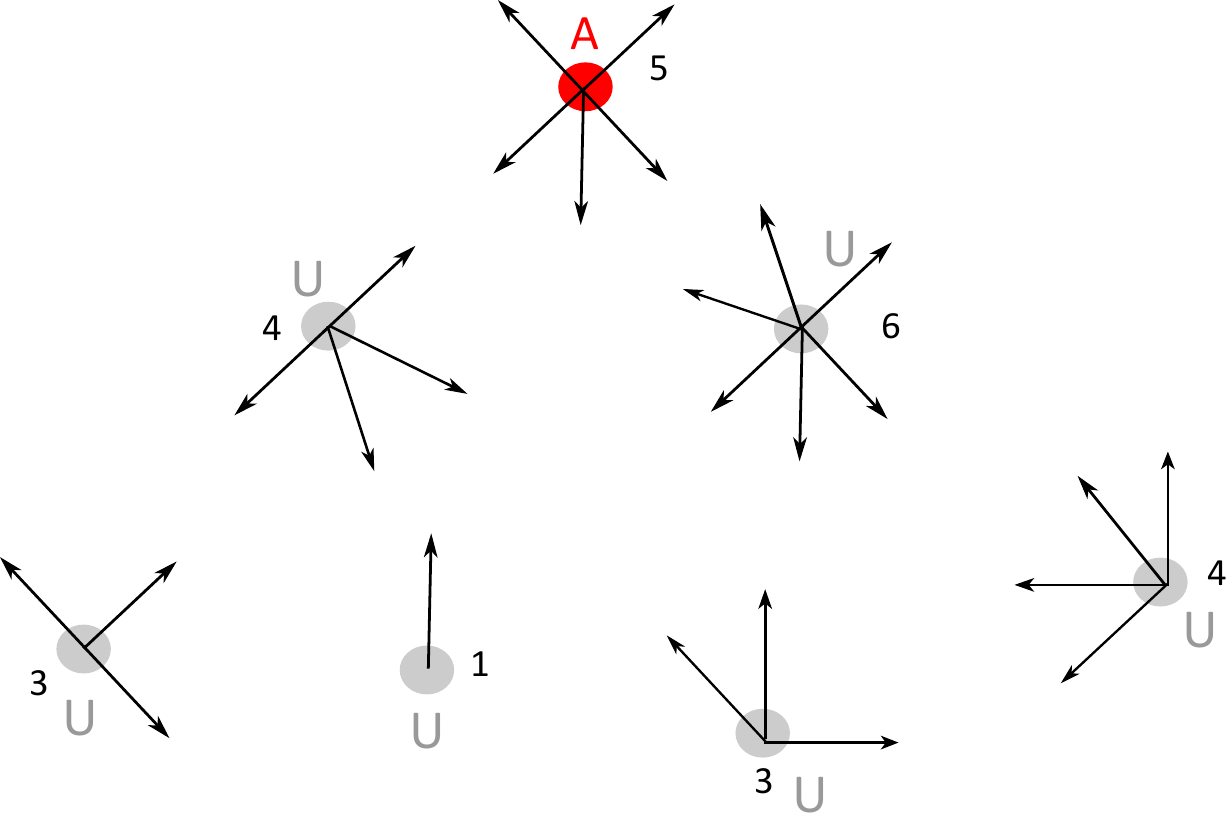}
\caption{Step 1 - the new active node is selected.}
\label{fig:exampleStep1}
\end{figure}
\begin{example}
\label{ex:construct}
Let the measure at $t^-$ be  
$$\mu_{t^-}=\delta_1+2\delta_3+2\delta_4+\delta_5+\delta_6,$$
so that the associated graph has $n=\cro{\mu_{t^-}, \mathbf 1}=7$ nodes and $\cro{\mu_{t^-},\chi}=26$ half-edges 
(see Fig. \ref{fig:exampleStep1}). A clock rings at $t$. 
The new $\sA$-node $a$ has degree $K\left(\mu_{t^-}\right)=5$ and the measure is updated to
$$\mu_{t^{-+}}=\delta_1+2\delta_3+2\delta_4+\delta_6.$$
\end{example}

\item[{\bf Step 2.}] \textbf{The neighbors of $a$ become of class} $\sB$: the neighbors of the new $\sA$-node are blocked, we say that they become 
$\sB$-nodes. 
The identity of these new $\sB$-nodes is determined by matching the $K\left(\mu_{t^-}\right)$ elements of $\maH_{t^-}$ of ego $a$, 
with half-edges of $\maH_{t^-}$, as follows:
\begin{itemize}
\item a first half-edge of $\maH_{t^-}$ of ego $a$ is matched with another one, drawn uniformly among the $\cro{\mu_{t^-},\chi}-1$ possible ones; 
\item on and on, as long as all half-edges of $\maH_{t^-}$ of ego $a$ have not been matched, 
we take one of those, and match it with another half-edge of $\maH_{t^-}$ which has not yet been matched, drawn uniformly in the latter set. 
\end{itemize}
Notice that at each step, we may match couples of half-edges emanating from $a$ together - hence creating self-loops around $a$. 
At the end of this procedure, we let $\tilde K\left(\mu_{t^-}\right)$ be the number of edges linking $a$ 
to other nodes. Clearly, $\tilde K\left(\mu_{t^-}\right)$ cannot exceed the number of half-edges of ego $a$, nor the number of half-edges of 
$\maH_{t^-}$ of egos different from $a$, in other words
\begin{equation}
\label{eq:boundtildeK}
\tilde K\left(\mu_{t^-}\right) \le K\left(\mu_{t^-}\right)\wedge \Bigl(\cro{\mu_{t^-},\chi}-K\left(\mu_{t^-}\right)\Bigl).
\end{equation} 
We have thus fixed the identity of the $q$ new $\sB$-nodes (where $q \le \tilde K\left(\mu_{t^-}\right)$), which are the 
egos $u_{t^{-}}(i_1),u_{t^{-}}(i_2),...,u_{t^{-}}(i_q)$ different from $a$, of the $\tilde K\left(\mu_{t^-}\right)$ half-edges matched with the $\tilde K\left(\mu_{t^-}\right)$ half-edges of ego $a$. 
We then set 
\[\left\{
\begin{array}{ll}
\maU_{t^{-++}}&=\maU_{t^{-+}} \setminus \left\{u_{t^{-}}(i_1),u_{t^{-}}(i_2),...,u_{t^{-}}(i_q)\right\};\\
\maA_{t^{-++}}&=\maA_{t^{-+}};\\
\maB_{t^{-++}}&=\maB_{t^{-+}}\cup\left\{u_{t^{-}}(i_1),u_{t^{-}}(i_2),...,u_{t^{-}}(i_q)\right\}.
\end{array}
\right.\] 
For all $j\in \llbracket 1,q \rrbracket$, let $N_{j}\left(\mu_{t^-}\right)$ be the number of edges shared by $u_{t^{-}}(i_j)$ with $a$. Let us  
define the two following point measures, 
\begin{align}
Y\left(\mu_{t^-}\right) &= \sum_{j=1}^q N_{j}\left(\mu_{t^-}\right)\delta_{d_{u_{t^{-}}(i_j)}};
\label{eq:defY}\\
\tilde Y\left(\mu_{t^-}\right) &=\sum_{j=1}^q \delta_{d_{u_{t^{-}}(i_j)}}.
\label{eq:defYtilde}
\end{align}  
In other words, for any $i$, $\tilde Y\left(\mu_{t^{-}}\right)(i)$ (resp., $Y\left(\mu_{t^{-}}\right)(i)$) is the number of neighbors of $a$ (resp., of edges shared by $a$ with its neighbors) having $i$ unmatched half-edges at $t^-$. 
Thus, $\cro{Y(\mu_{t^-}), \mathbf 1}=\tilde K(\mu_{t^-})$ is the number of half-edges of ego $a$ and $\cro{\tilde Y(\mu_{t^-}), \mathbf 1}=q$ is the number of neighbors of $a$.
As the new $\sB$-nodes are no longer unexplored, their degree must be erased from the measure $\mu_{t^{-+}}$,
which is updated as follows,   
\begin{equation}
\label{eq:genetape2}
\mu_{t^{-++}}=\mu_{t^{-+}}-\tilde Y\left(\mu_{t^-}\right).
\end{equation}
There remain $\cro{\mu_{t^-},\chi}-K\left(\mu_{t^-}\right)-\tilde K\left(\mu_{t^-}\right)$ unmatched half-edges at this point. 

\begin{figure}
\centering
\includegraphics[scale=0.9]{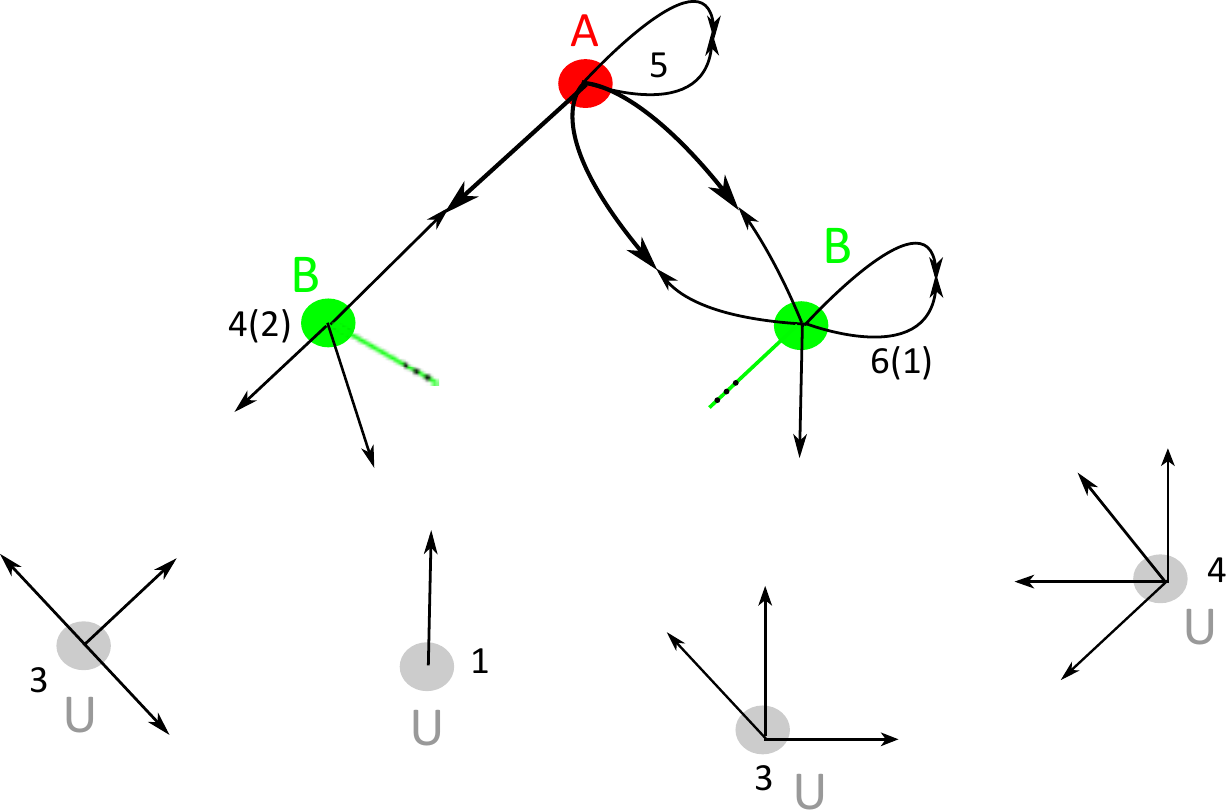}
\caption{Step 2 - the blocked nodes are attached to the active node.}
\label{fig:exampleStep2}
\end{figure}

\begin{example}[Example \ref{ex:construct} continued] 
The uniform selection of the neighbors of the new $\sA$-node $a$ results in a loop around it, 
so we have $\tilde K\left(\mu_{t^-}\right)=3$, $q=2$,
$$Y\left(\mu_{t^-}\right)=\delta_4 + 2\delta_6\quad \mbox{and} \quad \tilde Y\left(\mu_{t^-}\right)=\delta_4 + \delta_6,$$
as there exists a double-edge between $a$ and its neighbor of degree 6 (see Fig. \ref{fig:exampleStep2}).  
We then have 
$$\mu_{t^{-++}}=\delta_1+2\delta_3+\delta_4.$$
Observe that $\tilde K\left(\mu_{t^-}\right) = 3$ and the number of unmatched half-edges at this point is
$\cro{\mu_{t^-},\chi}-K\left(\mu_{t^-}\right)-\tilde K\left(\mu_{t^-}\right)= 26- 5 - 3=18.$
Between parenthesis is indicated the number of still unmatched half-edges of the $\sB$-nodes. 
\end{example}

\item[{\bf Step 3.}] \textbf{ Updating of the number of unmatched half-edges}: 
the available half-edges at this point, {\em i.e.} the elements of $\maH_{t^{-++}}$, either emanate from $\sB$-nodes and do not point to $a$, 
or emanate from nodes of $\maU_{t^{-++}}$. Let us denote 
$$\maH^{\scriptsize{\sB}}_{t^{-++}}:= \Bigl\{\mbox{half-edges of $\maH_{t^{-++}}$ having ego in $\maB_{t^{-++}}$}\Bigl\} \subset \maH_{t^{-++}},$$
and observe that
\begin{equation}
\label{eq:boundX} 
\left |\maH^{\scriptsize{\sB}}_{t^{-++}}\right |=\cro{\tilde Y\left(\mu_{t^-}\right),\chi} - \tilde K\left(\mu_{t^-}\right)
=\cro{\tilde Y\left(\mu_{t^-}\right),\chi} - \cro{Y\left(\mu_{t^-}\right),\mathbf 1}.
\end{equation}
We now fully attach the new $\sB$-nodes to the associated graph, {\em i.e.} we match all the half-edges of 
$\maH^{\scriptsize{\sB}}_{t^{-++}}$, either with other elements of $\maH^{\scriptsize{\sB}}_{t^{-++}}$, 
or with elements of $\maH_{t^{-++}}$ having ego in $\maU_{t^{-++}}.$ This is done according to the following procedure: 
\begin{itemize}
\item draw an integer uniformly at random in $\llbracket 1,q \rrbracket$, say $\ell$. 
We match the remaining $d_{u_{t^{-}}(i_\ell)}-N_{\ell}\left(\mu_{t^-}\right)$ open half-edges of $u_{t^{-}}(i_\ell)$ exactly as those of $a$: 
we take these open half-edges one by one; each time, we draw uniformly at random a match for the latter in all available half-edges (emanating from 
a node of $\maU_{t^{-++}}$, another node of $\maB_{t^{-++}}$ or from $u_{t^{-}}(i_\ell)$ itself), until all half-edges of $u_{t^{-}}(i_\ell)$ are matched;
\item then, draw at random another integer $m$ in $\llbracket 1,q \rrbracket\setminus \{\ell\}$, and match all available half-edges of 
$u_{t^{-}}(i_m)$ in the same manner, and so on... until all half-edges of $\maH^{\scriptsize{\sB}}_{t^{-++}}$ have been matched, to form edges 
of the associated graph. 
\end{itemize}
At the end of this operation, we have possibly created edges between the new $\sB$-nodes and the remaining $\sU$-nodes. 
Let us denote $X\left(\mu_{t^-}\right)$ the number of such edges, in other words 
\begin{equation}
\label{eq:defX}
X\left(\mu_{t^-}\right)=\mbox{Card}\,\Bigl\{\mbox{half-edges of $\maH^{\scriptsize{\sB}}_{t^{-++}}$ matched with half-edges of } 
\maH_{t^{-++}}\setminus \maH^{\scriptsize{\sB}}_{t^{-++}}\Bigl\}.
\end{equation}
To update $\mu_{t^{-++}}$, we have to subtract these $X\left(\mu_{t^-}\right)$ half-edges to the number of available half-edges of the remaining $\sU$-nodes. 
To formalize this operation, it is convenient to index the remaining $\sU$-nodes in the following way: for all $i\in \llbracket 1,n-1 \rrbracket$ 
and all $\ell \in \llbracket 1,\mu_{t^{-++}}(i) \rrbracket$, we let $u_{t^{-++}}(i,\ell)$ be the $\ell$-th node of $\maU_{t^{-++}}$ having $i$ unmatched half-edges at $t^{-++}$ (if any), ranked in an arbitrary order. 
Then, we define 
\begin{equation}
\label{eq:defW}
W\left(\mu_{t^-}\right)(i,\ell)=\mbox{Card}\,\Bigl\{\mbox{edges shared by $u_{t^{-++}}(i,\ell)$ with nodes of $\maB_{t^{-++}}$}\Bigl\},
\end{equation} 
and observe that 
$$\sum_{i=1}^{n-1}\sum_{\ell=1}^{\mu_{t^{-++}}(i)} W\left(\mu_{t^-}\right)(i,\ell)=X\left(\mu_{t^-}\right).$$
The quantity $W\left(\mu_{t^-}\right)(i,\ell)$ has to be subtracted from the number of open half-edges of each unexplored node $u_{t^{-++}}(i,\ell)$, hence we finally write   
\begin{equation}
\label{eq:genetape3}
\displaystyle\mu_{t}=\mu_{t^{-++}} - \sum_{i=1}^{n-1}\sum_{\ell=1}^{\mu_{t^{-++}}(i)}\left(\delta_{i}-\delta_{i-W\left(\mu_{t^-}\right)(i,\ell)}\right).
\end{equation}
\end{itemize}

\begin{figure}
\centering
\includegraphics[scale=0.9]{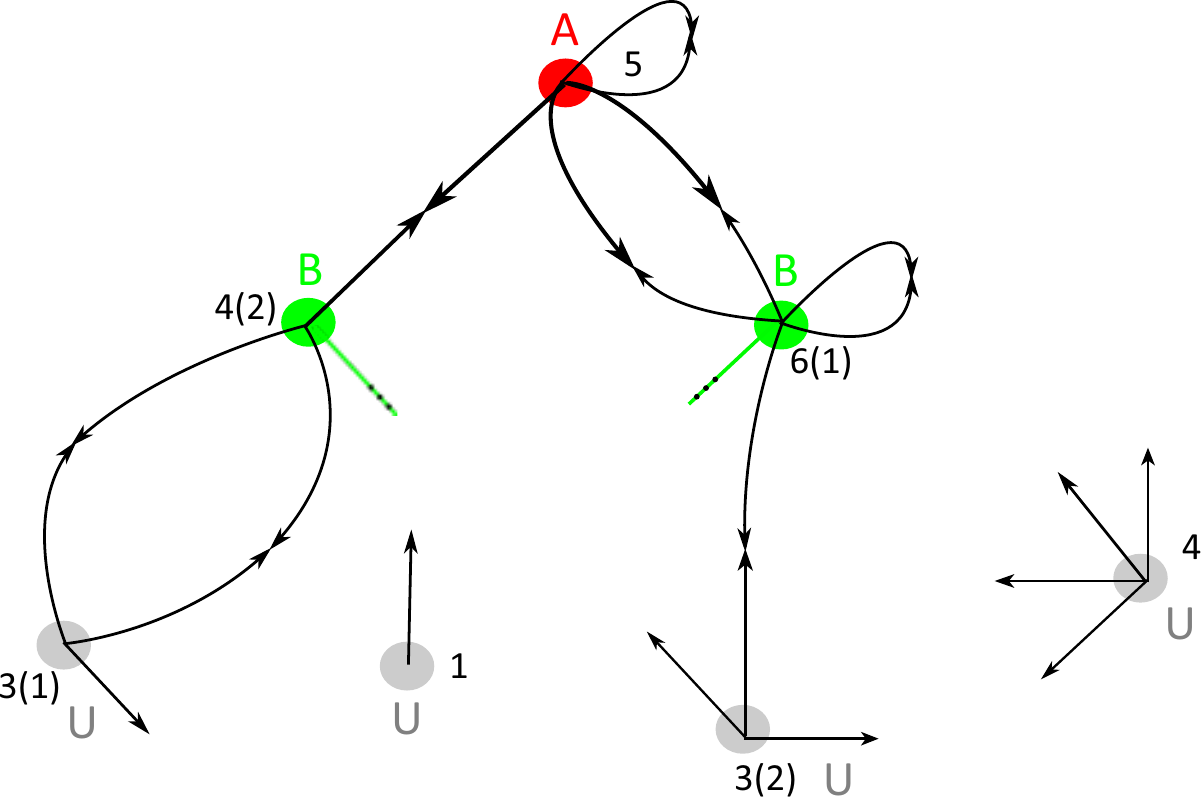}
\caption{Step 3 - the new blocked nodes are connected between each other, and with remaining unexplored nodes}
\label{fig:exampleStep3}
\end{figure}

\begin{example}[Example \ref{ex:construct} concluded]
In Fig. \ref{fig:exampleStep3} we see that $\left|\maH^{\scriptsize{\sB}}_{t^{-++}}\right|=\cro{\tilde Y\left(\mu_{t^+}\right),\chi}-\tilde K\left(\mu_{t^-}\right) = 10-3=7$ and $X\left(\mu_{t^-}\right)=3$, since there are 4 half-edges emanating from $\sB$-nodes that are matched together. 
The three remaining half-edges are matched with the remaining $\sU$-nodes as follows: if among the $\sU$-nodes of degree 3, the one on the left has label 1 and that on the right has label 2, then
$$W\left(\mu_{t^-}\right)(3,1)=2\mbox{ and }W\left(\mu_{t^-}\right)(3,2)=1.$$
The updated measure is then 
$$\mu_{t}=2\delta_1+\delta_2+\delta_4.$$
There remain $\cro{\mu_t, \mathbf 1}=4$ unexplored nodes and $\cro{\mu_t,\chi}=8$  unmatched half-edges. 
\end{example}

After Steps 1-3, we end up with a measure $\mu_{t}$ where $1+q$ atoms have been erased with respect to 
$\mu_{t^-}$, and whose first moment ({\em i.e.} the number of open half-edges) has the same parity as that of $\mu_{t^-}$. 
The associated graph $G_t$ equals $G_{t^-}$, plus all edges that have been drawn between $a$ and its neighbors, between its neighbors with one another, 
and between its neighbors and the remaining $\sU$-nodes. Moreover, only the remaining $\sU$-nodes still have open half-edges at $t$ and the measure 
$\mu_t$ provides the repartition of the latter half-edges among $\maU_t$. At this point, we re-index all elements of $\maU_t$ in the order 
of increasing number of open half-edges, as was done above: 
$$\maU_{t}=\left\{u_{t}(1),...,u_{t}(U_t)\right\}.$$
We now draw a new exponential r.v. $\xi_t$ of parameter $\lambda U_t$, independently of everything else. 
As above, the system remains constant until time $t+\xi_t$, at which we re-iterate Step 1-3, and so on.

The procedure ends at the stopping time 
$$T_0=\inf\Bigl\{t\ge 0;\,U_t=0\Bigl\},$$
which clearly is a.s. finite. At that instant, if $\sum_{i=1}^n d(i)$ was odd there remain a single unmatched 
half-edge, which we remove. At that time $T_0$, we thus end up with $\mu_{T_0}=\mathbf 0$, $\maU_{T_0}=\maH_{T_0}=\emptyset$ and 
$\left|\maA_{T_0}\cup\maB_{T_0} \right|=n$. The final associated graph $G_{T_0}$ is a multi-graph of degree vector $\mathbf d$ (up to the deletion of a single half-edge in the case mentioned above). The set $\maA_{T_0}$ is the {\em jamming limit} of the latter graph (and of the parking process) 
and $\left|\maA_{T_0}\right|/n$, its {\em jamming constant}.

\section{Generator and semi-martingale decomposition}\label{sec:generator}

It follows from the construction of Section \ref{sec:construct}, that the process $\proc{\mu}$ is Markov in $\M_F(\N)$. 
Its infinitesimal generator $\mathscr Q$ is defined for all $F:\M_F(\N) \to \mathbb R $ in the domain of $\mathscr Q$, by 
\begin{gather}
\label{eq:defgenerateur}
\mathscr Q F(\mu)=\lim_{h \to 0} \frac{1}{h}\biggl(\esp{F(\mu_h) \mid \mu_0=\mu}-F(\mu)\biggl),\,\mu \in \M_F(\N).
\end{gather}

Let $\phi:\R\to \R$ be a bounded function, and denote for all $\mu \in \M_F(\N)$, $\Pi_\phi(\mu)=\cro{\mu,\phi}.$
We show that $\Pi_\phi$ belongs to the domain of $\mathscr Q$ and deduce $\mathscr Q \Pi_\phi$ from (\ref{eq:genetape1},\ref{eq:genetape2},\ref{eq:genetape3}), 
following Steps 1-3.
 
{\bf Step 1.} The probability that a node of degree $k$ is drawn, equals the proportion of atoms at level $k$ among all atoms of the measure $\mu_0$. In other words, for any $\mu$ and $k$, 
\begin{equation}
\label{eq:genetape1bis}
\P_K(k \mid \mu):=\pr{K(\mu_0)=k \mid \mu_0=\mu}=\frac{\mu(k)}{\cro{\mu,\mathbf 1}}.
\end{equation}
Denote for any $\mu$, $k$, and any integer $\tilde k \le k\wedge(\cro{\mu,\chi}-k)$, 
\begin{equation}
\label{eq:distribtildek}
\P_{\tilde K}(\tilde k \mid \mu,k):=\pr{\tilde K(\mu_0)=\tilde k \mid \mu_0=\mu,\,K(\mu)=k}
\end{equation}
and $\mathbb E_{\tilde K}[. \mid \mu,k]$ the corresponding conditional expectation. 
We won't need a precise expression of the latter distribution, but will show instead that 
$\tilde K$ coincides with $K$ with overwhelming probability as the size $n$ of the graph goes large. Or which is equivalent, when the number of nodes goes to infinity the probability of choosing two half-edges of the new $\sA$-node to match together is arbitrarly small.

To see this, consider a system in which all half-edges emanating from the new $\sA$-node $a$ are duplicated, and where we match the $K(\mu_0)$ half-edges of $a$ uniformly, either with half-edges emanating from the other nodes, or with their duplicata. Hence, given $\mu$ and $k$, in this auxiliary system the number $\breve K(\mu_0)$ of half-edges of ego $a$ matched with half-edges of other egos (and not with their duplicata) follows an hypergeometric distribution with parameters $k$, $\cro{\mu,\chi}$ and $\cro{\mu,\chi}-k$. 
Moreover the probability of matching as many half-edges of ego $a$ with half-edges of other egos as possible 
(\emph{i.e.} $K(\mu_0) \wedge \left(\cro{\mu_0,\chi}-K(\mu_0)\right)$ half-edges) is clearly larger in the actual system than in the auxiliary one, thus 
\begin{multline}
\mathbf P\Bigl[\mbox{the new $\sA$-node at $0$ shares its
$K(\mu_0)$ edges with other nodes}\mid \mu_0=\mu,\,K(\mu)=k \Bigl]\\
\begin{aligned}
        &=\P_{\tilde K}\left(k \mid \mu,k\right)\\
        &\ge \pr{\breve K(\mu_0) = k  \mid \mu_0=\mu,\,K(\mu)=k}\\
        &={{\cro{\mu,\chi}-k  \choose k} \over {\cro{\mu,\chi} \choose k}}\ind_{\{k \le \cro{\mu,\chi}-k\}}, 
\end{aligned}
\label{eq:Pkktilde}
\end{multline}
which will be proven to tend to 1 as $n$ goes large.

{\bf Step 2.} Conditionally to $\{\mu_0=\mu\} \cap \{K\left(\mu_0\right)=k\}\cap \{\tilde K\left(\mu_0\right)=\tilde k\}$,
$Y\left(\mu_0\right)$ follows a multivariate hypergeometrical distribution on $\M_F\left(\N\right)$, of parameters
$\left(\tilde k,\cro{\mu,\chi}-k,n-1,P\right)$ (see \ref{sec:combi}), where $P$ is given 
by
$$P(i)=i \left(\mu (i) - \delta_k(i)\right),\,i\in\llbracket 0,n-1 \rrbracket.$$
In other words, we have for all $y \in \M_F\left(\llbracket 0,n-1 \rrbracket\right)$ such that $\cro{y,\mathbf 1} = \tilde k$ and
$y(i)\le P(i)$ for all $i$,

\begin{align}
\P_Y(y \mid \mu,k,\tilde k) &:= \pr{Y(\mu_0)=y \mid \mu_0=\mu,\,K\left(\mu_0\right)=k,\,\tilde K\left(\mu_0\right)=\tilde k}\nonumber\\
&={\displaystyle\prod_{i\in\llbracket 0,n-1 \rrbracket} {i\left(\mu(i)-\delta_{k}(i)\right) \choose y(i)}\over {\cro{\mu,\chi}-k \choose \tilde k}},
\label{eq:distriby}
\end{align}
and we denote $\mathbb E_Y\left[ . \mid \mu,k,\tilde k\right]$ the corresponding conditional expectation.
For all $\tilde y \in \M_F(\N)$, we denote
\begin{equation}
\label{eq:distribtildey}
\P_{\tilde Y}(\tilde y \mid \mu,k,\tilde k,y):=\pr{\tilde Y(\mu_0)=\tilde y \mid \mu_0=\mu,\,K(\mu)=k, \tilde K(\mu)=\tilde k,\,Y(\mu)=y},
\end{equation}
and define $\mathbb E_{\tilde Y}\left[ . \mid \mu,k,\tilde k,y\right]$ accordingly. Here again, the precise form of this distribution is not needed explicitly, since $\tilde{Y}(.)$ and $Y(.)$ tend to coincide for large graphs, {\em i.e.} $\P_{\tilde Y}(y \mid \mu,k,\tilde k,y)$ (the probability that all edges between the new active node and its neighbors are simple) tends to 1 for large $n$.
This will be shown using the following lower-bound: recalling (\ref{eq:notationv}), we have that
\begin{align}
\P_{\tilde Y}\left(y \mid \mu,k,\tilde k,y\right)
&=\sum_{\substack{\beta \in \left(\{1,...,\cro{\mu,\mathbf 1}-1\}\right)^{\tilde k};\\\beta(1)\ne \beta(2)\ne ...\ne \beta(k)}}\frac{\prod_{i=1}^{\tilde k}
                                  d_{v_{\beta(i)}}\left(\mu\right)}{\maA^{\tilde k}_{\cro{\mu,\chi}-k}}\nonumber\\
&=\frac{\prod_{i\in \N^*}i^{y(i)}}{\maA^{\tilde k}_{\cro{\mu,\chi}-k}}\nonumber\\
&\ge 1-{\tilde k \choose 2}
     \sum_{\ell =1}^{ \cro{\mu,\mathbf 1}-1}\frac{ d_{v_\ell}\left(\mu-\delta_k\right)\left(d_{v_\ell}\left(\mu-\delta_k\right)-1\right)}
{\left(\cro{\mu,\chi}-k\right)\left(\cro{\mu,\chi}-k-1\right)}\nonumber\\
&\ge 1-{\tilde k \choose 2}\frac{\cro{\mu,\chi^2}}{\left(\cro{\mu,\chi}-k\right)\left(\cro{\mu,\chi}-(k+1)\right)},
\label{eq:majorePytilde}
\end{align}
since the number of (ordered) configurations entailing a multiple edge between $a$ and a given $\sU$-vertex $u$,
is less than the number of pairs of half edges starting from $u$, times the number of couples of half-edges matched with the $\tilde k$ half-edges of $a$.

{\bf Step 3.} Let us first address the r.v. $X(.)$ defined by (\ref{eq:defX}).  
Denote for all $y,\tilde y$ and all integers $x \le \cro{\tilde y,\chi}-\cro{y,\mathbf 1}$,
\begin{equation}
\label{eq:distribX}
\P_X(x \mid \mu,k,\tilde k,y,\tilde y):=\pr{X\left(\mu_0\right)=x \mid \mu_0=\mu,\,K\left(\mu\right)=k,\,\tilde K\left(\mu\right)= \tilde k
,\,Y\left(\mu\right)=y,\,\tilde Y\left(\mu\right)=\tilde y}.
\end{equation}
Again, the exact distribution is not needed, as we will see that $X(\mu_0)$ tends to coincide with the number of open half-edges emanating from the new 
$\sB$-nodes (in other words, the latter nodes tend not to share any edge with one another - which would have created triangles including $a$). 
Indeed, noticing that the total number of available half-edges after the connexion between the new $\sA$ node and its neighbors equals 
$\cro{\mu_0,\chi}-K\left(\mu_0\right)-\tilde K\left(\mu_0\right)$ and using (\ref{eq:boundX}), 
we can show exactly as for (\ref{eq:Pkktilde}) that 
\begin{align}
\P_X(\cro{\tilde y,\chi}-\tilde k \mid \mu,k,\tilde k,y,\tilde y)&=\mathbf P\Bigl[\mbox{the new $\sB$-nodes do not share any edge with one another}\nonumber\\
 &\quad\quad\quad\quad\quad\quad\quad  \mid \mu_0=\mu,\,K(\mu)=k,\,\tilde K(\mu)=\tilde k,\,Y(\mu)=y,\,\tilde Y(\mu)=\tilde y\Bigl]\nonumber\\
        &\ge {{\cro{\mu,\chi}-k - \cro{\tilde y,\chi} \choose \cro{\tilde y,\chi}-\tilde k} \over 
               {\cro{\mu,\chi} - k - \tilde k  \choose \cro{\tilde y,\chi}-\tilde k}}\ind_{\left\{\cro{\tilde y,\chi}-\tilde k \le  \cro{\mu,\chi}-k - \cro{\tilde y,\chi}\right\}},\label{eq:QX}
\end{align}
which is the probability that an hypergeometric random variable with parameters $\cro{\tilde y,\chi}-\tilde k$,  $\cro{\mu,\chi}-k-\cro{\tilde y,\chi}$ and 
$\cro{\mu,\chi}-k-\tilde k$  takes the value $\cro{\tilde y,\chi}-\tilde k$.
As above, the latter quantity will be shown to tend to one for large values of $n$.

Now, let $\mu,k,\tilde k,y,\tilde y$ and $x$ as above, and define the set 
\begin{multline*}
\mathcal W(\mu,k,\tilde k,y,\tilde y,x):=\Biggl\{w \in \N^{\llbracket 1,n-1 \rrbracket^2};\sum_{i=1}^{n-1}\sum_{\ell=1}^{n-1}w_{i\ell}=x
\mbox{ and for all }
i\in\llbracket 1,n-1 \rrbracket,\\w(i,\ell) \le i\mbox{ and }w(i,\ell)=0\mbox{ for all }\ell > \mu(i)-\delta_k(i)-\tilde y_i.\Biggl\}.
\end{multline*}
Clearly, the double-indexed sequence $W(\mu)$ defined by (\ref{eq:defW}) is drawn according to an hypergeometric choice among the ``bunches" of half-edges 
represented by the remaining $\sU$-nodes. So we obtain readily that for all $w \in \mathcal W(\mu,k,\tilde k,y,\tilde y,x)$, 
\begin{multline}
\P_W(w \mid \mu,k,\tilde k, y, \tilde y,x)\\
=\pr{W\left(\mu_0\right)=w \mid \mu_0=\mu,\,K\left(\mu_0\right)=k,\,\tilde K\left(\mu_0\right)=\tilde k,\,Y\left(\mu_0\right)=y,\,\tilde Y\left(\mu_0\right)=\tilde y,\,X\left(\mu_0\right)=x}\\
=\displaystyle\frac{\prod_{i\in\llbracket 1,n-1 \rrbracket}\prod_{\ell =1}^{\mu(i) - \delta_k(i)-\tilde y(i)}{i \choose w_{i\ell}}}{{\cro{\mu,\chi}-k-\cro{\tilde y,\chi} \choose x}}.\label{eq:genetape3ter}
\end{multline}
We also need to introduce the following probability,
\begin{multline}
\label{eq:defcheckQ}
\mathbb Q_W\left(\mu,k,\tilde k,y,\tilde y,x\right)=\sum_{w;w(i,\ell)\le 1\forall i,\ell}\P_W(w \mid \mu, k, x, y)\\
\shoveleft{=\mathbf P\Biggl[\mbox{no new \sB-vertex at 0 has a multiple edge towards a remaining $\sU$-node}}\\
                                          \biggl| \mu^n_0=\mu,\,K\left(\mu_0\right)=k,\,\tilde K\left(\mu_0\right)=\tilde k
,\,Y\left(\mu_0\right)=y\,\tilde Y\left(\mu_0\right)=\tilde y,\,X\left(\mu^n_0\right)=x\Biggl].
\end{multline}
Just as (\ref{eq:majorePytilde}), we obtain that
\begin{multline}
\mathbb Q_W\left(\mu,k,\tilde k,y,\tilde y,x\right)
\ge 1-{x \choose 2}\sum_{\ell =1}^{ \cro{\mu,\mathbf 1}-1-k}\frac{ d_{v_\ell}\Bigl(\mu-\delta_k-\sum_i \tilde y(i)\delta_i \Bigl)\left(d_{v_\ell}\Bigl(\mu-\delta_k-\sum_i \tilde y(i)\delta_i\Bigl)-1\right)}{\left(\cro{\mu,\chi}-k-\cro{y,\chi}\right)\left(\cro{\mu,\chi}-k-\cro{y,\chi}-1\right)}\\
\ge 1-{x \choose 2}\frac{ \cro{\mu,\chi^2}}{\left(\cro{\mu,\chi}-k-\cro{y,\chi}\right)\left(\cro{\mu,\chi}-k-\cro{y,\chi}-1\right)}
.\label{eq:majorecheckQ}
\end{multline}
We finally introduce, for all $i\in\llbracket 1,n-1 \rrbracket$, 
\begin{equation}
\label{eq:defZ}
Z\left(\mu_0\right)(i)=\sum_{\ell=1}^{n-1}W(i,\ell)\left(\mu_0\right),
\end{equation}
the number of half-edges of $\sU$-nodes having initially $i$ unmatched half-edges, that have been matched with half-edges emanating from the new $\sB$-nodes. 
It then follows from (\ref{eq:genetape3ter}) that, given $\mu,k,\tilde k,y,\tilde y$ and $x$, 
the integer measure $Z\left(\mu_0\right)$ follows a multivariate hypergeometric distribution of parameters
$\left(x,\cro{\mu,\chi}-k-\cro{\tilde y,\chi},n-1,P^\prime\right)$, where
$$P^\prime(i)=i\left(\mu(i)-\delta_k(i)- \tilde y(i)\right) ,\,i\in \llbracket 1,n-1 \rrbracket.$$
In other words, for all $z\in\M_F(\N)$,
\begin{multline}
\label{eq:distribz}
\P_W\left(Z(\mu)=z|\mu,k,\tilde k, y,\tilde y, x\right)\\
=\pr{Z(\mu_0)=z \mid \mu_0=\mu,\,K\left(\mu_0\right)=k,\,\tilde K\left(\mu_0\right)=\tilde k,\,Y\left(\mu_0\right)=y,\,
\tilde Y\left(\mu_0\right)=\tilde y,\,X\left(\mu_0\right)=x}\\
=\frac{\displaystyle\prod_{i\in\llbracket 1,n-1 \rrbracket} {i\left(\mu(i)-\delta_k(i)-\tilde y(i)\right) \choose z(i)}}{{\cro{\mu,\chi}-k-\cro{\tilde y,\chi} \choose x}}\cdot
\end{multline}

\noindent Therefore, for all $h>0$, from (\ref{eq:genetape1},\ref{eq:genetape2},\ref{eq:genetape3}) and (\ref{eq:genetape1bis}) we obtain that for all $\mu\in\M_F(\N)$,
\begin{multline}
\mathscr Q \Pi_{\phi}(\mu)=\lim_{h\to 0}\frac{1}{h}\esp{\left(\cro{\mu_h,\phi}-\cro{\mu,\phi}\right)\ind_{\{\xi_0<h\}} \,\mid\,\mu_0=\mu}\\
\shoveleft{=\lambda\cro{\mu,\mathbf 1}\Biggl [-\sum_{k\in\N}\frac{\mu(k)}{\cro{\mu,\mathbf 1}}}\\
\shoveleft{\times\left\{\phi(k)+ \sum_{\tilde k\in \llbracket 0,k\rrbracket} \P_{\tilde K} (\tilde k| \mu,k) \sum_{\substack{y \in \M_F(\N);\\
\cro{y,\mathbf 1} = \tilde k}}\P_Y(y \mid \mu,k,\tilde k)\sum_{\tilde y \in \M_F(\N)}\P_{\tilde Y}(\tilde y \mid \mu,k,\tilde k,y)\right.}\\
\shoveleft{\times\sum_{i=1}^{n-1}\Biggl(\tilde y(i)\phi(i)}\\
\shoveright{+  \sum_{x \in \N^*}\P_X(x \mid \mu,k,\tilde k,y, \tilde y)\sum_{w\in\N^{\llbracket 1,n-1 \rrbracket^2}}\P_W(w \mid \mu,k,\tilde k,y,\tilde y,x)\sum_{\ell=1}^{n-1}
\left(\phi(i)-\phi(i-w(i,\ell))\right)\Biggl)\Biggl\}\Biggl ]}\\
\shoveleft{=\lambda\Biggl [-\cro{\mu,\phi}-\sum_{k\in\N}\mu(k)}\\
\shoveleft{\times\left\{\sum_{\tilde k\in \llbracket 0,k\rrbracket}\P_{\tilde K} (\tilde k| \mu,k)\sum_{\substack{y \in \M_F(\N);\\\cro{y,\mathbf 1} 
=\tilde k}}\P_Y(y \mid \mu,k,\tilde k)\sum_{\tilde y \in \M_F(\N)}\P_{\tilde Y}(\tilde y \mid \mu,k,\tilde k,y)\right.}\\
\shoveleft{\times\sum_{i=1}^{n-1}\Biggl(\tilde y(i)\phi(i)}\\
+\sum_{x \in \N^*}\P_X(x \mid \mu,k,\tilde k,y,\tilde y)\sum_{w\in\N^{\llbracket 1,n-1 \rrbracket^2}}\P_W(w \mid \mu,k,\tilde k,y,\tilde y,x)\sum_{\ell=1}^{n-1}
\left(\phi(i)-\phi(i-w(i,\ell))\right)\Biggl)\Biggl\}\Biggl ]
%
\label{eq:gener1},
\end{multline}
where $\P_{\tilde K} (\tilde k| \mu,k)$, $\P_Y(y \mid \mu,k,\tilde k)$, $\P_{\tilde Y}(\tilde y \mid \mu,k,\tilde k,y)$, $\P_X(x \mid \mu,k,\tilde k,y,\tilde y)$ and $\P_W(w \mid \mu,k,\tilde k,y,\tilde y,x)$ are respectively defined by \eqref{eq:distribtildek}, \eqref{eq:distriby}, \eqref{eq:distribtildey}, \eqref{eq:distribX}
 and \eqref{eq:genetape3ter}.

\subsubsection*{Semi-martingale decomposition}

 Let $\proc{\mathcal F}$ be the natural filtration of $\proc{\mu}$. From (\ref{eq:gener1}) it is easy to check using Lemma 3.5.1 and Corollary 3.5.2 of 
\cite{crm} that $\proc{\mu}$ is a Feller-Dynkin process of the space $\mathcal D\left([0,\infty),\M_F(\N)\right)$.
It then follows from standard stochastic calculus that for any bounded $\phi: \R \to \R$, the following is a $\mathcal F$-martingale:


\begin{equation}
t\mapsto M(\phi)_t =\cro{\mu_t,\phi}-\cro{\mu_0,\phi} - \int_0^t \mathscr Q\Pi_\phi\left(\mu_s\right) \,ds,\label{eq:mart}
\end{equation}
of quadratic variation process given by
\begin{multline}
\label{eq:crochet}
t\mapsto\langle\!\langle M(\phi) \rangle\!\rangle_t=\int_0^t \biggl(\mathscr Q \left(\Pi_\phi\right)^2(\mu_s) -2\Pi_\phi(\mu_s)
\mathscr Q\Pi_{\phi}\left(\mu_s\right)\biggl)\,ds\\
\shoveleft{=\lambda\int_0^t \left [\sum_{k\in\N}\mu_s(k)\left\{\phi(k)
+\sum_{\tilde k\in \llbracket 0,k\rrbracket}\P_{\tilde K} (\tilde k| \mu_s,k)\sum_{y \in \M_F(\N)}\P_Y(y \mid \mu_s,k)\right.\right.}\\
\times\sum_{i=1}^{n-1}\left(\sum_{\tilde y \in \M_F(\N)}\P_{\tilde Y}(\tilde y \mid \mu_s,k,y)\tilde y(i)\phi(i)\right.\\
+\left.\left.\sum_{x \in \N^*}\P_X(x \mid \mu_s,k,y)\sum_{w\in\N^{\llbracket 1,n-1 \rrbracket^2}}\P_W(w \mid \mu_s,k,x,y)\sum_{\ell=1}^{n-1}
\left(\phi(i)-\phi(i-w_{i,\ell})\right)\Biggl)\right\}^2\right]\,ds.
\end{multline}

\section{Main results and consequences}\label{sec:hyd}

\subsection{Hydrodynamic limit}

We are interested in the behavior of the measure-valued process $(\mu_t)_{t\geq 0}$ as the size of the graph grows to infinity. 
We consider a sequence of models, where the size of the $n$-th graph equals $n$, and add a superscript $^n$ to all 
parameters and processes relative to the $n$-th system. 
Then, we scale the $n$-th process of empirical degree distributions as follows: 
$$\bar\mu^{n}_t=\frac{1}{n} \mu^n_t,\,t\ge 0$$
and we denote the normalized versions of the derived processes accordingly, {\em i.e.} for all $t\ge 0$, 
$$\bar U^n_t={U^n_t \over n};\quad\bar H^n_t={H^n_t \over n};\quad \bar M^{n}(\phi)_t=\frac{1}{n} M^{n}(\phi)_t,\,\phi \in \mathcal B_b.$$
Our main result is the following.
\begin{theorem}\label{theo.main}
Assume that for all $\phi\in \mathcal B_b \cup \{\chi^6\}$,
\begin{equation}
\label{eq:convinit}
\cro{\bar\mu^n_0,\phi} \tendP \cro{\zeta,\phi},
\end{equation}
where $\zeta$ is a deterministic element of $\M_F(\N)$ such that 
\begin{equation}
\label{eq:condinit}
0<\cro{\zeta,\chi}\,\mbox{ and }\,\cro{\zeta,\chi^6} <\infty.
\end{equation}
Then, for all $T$ and all $\phi \in \maB_b$ we have
$$\sup_{t\in[0,T]}\left|\cro{\barmun_t,\phi}-\cro{\barmu_t,\phi}\right| \tendP 0,$$
where $\bar \mu$ is the unique element of $\C\left(\R+,\M_F(\N)\right)$
 satisfying the following infinite dimensional differential system: for all $t\ge 0$ and all bounded $\phi$,
\begin{align}
\cro{\barmu_0,\phi}&=\cro{\zeta,\phi}\nonumber;\\
{d \over dt} \cro{\bar\mu_t, \phi} &= \cro{\Psi\left(\bar\mu\right)_t,\phi}\nonumber \\
 &:=
\left\{\begin{array}{ll}
-\lambda \left[\cro{\bar\mu_t,\left(\mathbf 1+\chi\right)\phi}+\cro{\bar\mu_t,\chi\Delta\phi}\left({\cro{ \bar\mu_t, \chi^2} \over \cro{\bar\mu_t,\chi }}-1
\right) \right]&\mbox{ if }\cro{\barmu_t,\chi}>0;\\
-\lambda \barmu_t(0)\phi(0)&\mbox{ if }\cro{\barmu_t,\chi}=0.
\end{array}\right.
\label{eq:limflu}
\end{align}
\end{theorem}

\begin{remark}
A consequence of the proof of this theorem is the existence of a solution to (\ref{eq:limflu}).
Its uniqueness is proved separately.
\end{remark}

\begin{remark}
By our very assumptions, $\cro{\mu^n_0,\mathbf 1}=n$ and thus $\cro{\barmun_0,\mathbf 1}=1$ for all $n$. It thus follows from (\ref{eq:convinit}) and (\ref{eq:condinit}) that $\barmu_0$ is a probability measure. 
In particular we have $$\cro{\barmu_0,\mathbf 1}\vee \cro{\barmu_0,\chi^6}< \infty,$$ a fact that will be used at several points of the proofs.
\end{remark}

\subsection{Main characteristics, and the jamming constant}
\label{subsec:particular}
It follows from Theorem \ref{theo.main} that the sequence 
$\{\bar U^n\}$ tends uniformly over compact time sets to the deterministic functions $\bar u$  
given for all $t$ by  $\bar u_t=\cro{\bar\mu_t,\mathbf 1}$ and which, from (\ref{eq:limflu}), 
satisfies 
\begin{equation}
\label{eq:diffu}
\dot{\bar u}_t =-\lambda \left(\bar u_t + \bar h_t\right),\,t \ge 0.
\end{equation}
On another hand, as a simple consequence of 
Theorem \ref{theo.main}, we have the weak convergence 
$$\barmun \Rightarrow \barmu \mbox{ in }\mathcal D\left(\R+,\M_F(\N)\right).$$ 
It then follows from the continuity of the mapping 
\[\left\{\begin{array}{ccc}
 \D\left(\R+,\M_F(\N)\right) &\longrightarrow &\D\left(\R+,\R\right)\\
 \mu &\longmapsto &\cro{\mu_.,\chi}
             \end{array}
             \right.\]
and from the Continuous Mapping Theorem, that $\bar H^n \Rightarrow \bar h$ in $\D\left([0,T],\R\right)$, where 
$\bar h_t=\cro{\bar\mu_t,\chi}$ for all $t\ge 0$. Again from (\ref{eq:limflu}), we can easily check that 
$$
\dot{\bar h}_t =-2\lambda \cro{\bar \mu_t,\chi^2},\,t\ge 0.
$$

\paragraph{Jamming constant}
Denote for all $n \in \N^*$ and all $t\ge 0$,
\begin{equation}
\label{eq:defJn}
J^n_t:=\left|\maA^n_t \right|
\end{equation}
the number of active nodes at $t$. 
The {\em jamming constant $\bar J^n$} of the associated graph, is the proportion of active nodes at the ending time $T^n_0$ of the 
exploration process. 
In other words, it is given by
\begin{equation}
\label{eq:jammingconstant}
\bar J^n={J^n_{T^n_0} \over n}.
\end{equation}
The following result can be deduced from Theorem \ref{theo.main}. Its proof is provided in section \ref{sec:proofjamming}.
\begin{coro}[Jamming constant of random graphs]\label{prop:jamming}
Under the assumption of Theorem \ref{theo.main}, we have that 
$$\bar J^n \tend c_{\zeta}\mbox{ in }L^1,$$
where
\begin{equation}
\label{eq:cmu0}
c_{\zeta}= \lambda \int_0^\infty \bar u_t\,dt=\lambda\int_0^\infty \cro{\bar\mu_t,\mathbf 1} dt
\end{equation}
and $\proc{\barmu}$ is the only solution to (\ref{eq:limflu}). 
\end{coro}
\begin{remark}
As expected, it readily follows from (\ref{eq:cmu0}) and (\ref{eq:diffu}) that  $c_{\zeta}= \int_0^\infty \bar u_{t/\lambda}\,dt $ does not depend on
$\lambda$ - we hence fix $\lambda=1$ in section
\ref{subsec:jamming} without loss of generality.
\end{remark}

\subsection{Connexion with the configuration model}

As far as the structure of the associated random graph is concerned, our construction 
 mimics the so-called {\em configuration model} (or uniform model).
More precisely, denote $\text{{\sc CM}}(n,\mathbf d^n)$ the random (multi-)graph obtained by the uniform mapping of 
half-edges in a graph of $n$ nodes having degree vector $\mathbf d^n$, as described in \cite{remco, wormaldCM}. We have the following. 

\begin{prop}
\label{prop:equiv}
The associated random multi-graph constructed jointly with the exploration process in Section \ref{sec:construct} 
equals $\text{{\sc CM}}(n,\mathbf d^n)$ in distribution. 
\end{prop}

\bp
The result follows from the so-called {\em independence property} of the configuration model (see {\em e.g.} \cite{wormaldCM}): 
choosing whatever rule for matching the half-edges as long as a given half-edge is matched uniformly among all the unmatched half-edges at each step, provides a realization of $\text{{\sc CM}}(n,\mathbf d^n)$. This is exactly what is done here: at each step, first all half-edges of the new active node $a$ are matched 
with other half-edges that are chosen uniformly among available ones (including other half-edges of $a$), and then all the remaining open half-edges of the neighbors $b_1$, and then $b_2$, $b_3$,...and finally $b_q$ ($q$ being the number of neigbors of $a$) are matched according to the same rule. Therefore, at the end of the algorithm the associated multi-graph we have drawn is nothing but a realization of $\text{{\sc CM}}(n,\mathbf d^n)$. 
\ep

We can now link our construction with the more usual one, consisting in fixing the graph beforehand, and then 
building an independent set on the latter. 
The {\em parking process} on a uniform graph, mentioned in the introduction, can be formalized as follows: we first fix 
$\tilde G$, a realization of $\text{{\sc CM}}(n,\mathbf d^n)$. We then construct sequentially the independent set on $\tilde G$ according to the following procedure,  
\begin{itemize}
\item[(i)] At time 0, set all nodes of $\tilde G$ as unexplored, and initiate an exponential clock of intensity $\lambda n$;
\item[(ii)] Each time an exponential clock rings, select a new active node $a$ uniformly at random among all unexplored nodes;
 \item[(iii)] The neighbors of $a$ in $\tilde G$ all become of class $\sB$;
 \item[(iv)] We set another clock of intensity $\lambda U$, where $U$ is the cardinality of the set of unexplored nodes at that instant, and go to step (ii).
\end{itemize}
The algorithm is terminated as soon as the set of unexplored nodes of $\tilde G$ is empty, and the jamming constant of the graph is obtained as the  
proportion of active nodes at that time. 

Let for all $t$, $\tilde{\maU}^n_t$ the set of unexplored nodes at $t$ and for any $j \in \tilde{\maU}^n_t$, $d_j\left(\tilde{\maU}^n_t\right)$ the number of neighbors of $j$ in $\tilde{\maU}^n_t$ at $t$. Let also the following random point measure, 
\begin{equation}
\label{eq:defmutilde}
\tilde \mu^n_t=\sum_{j\in \tilde{\maU}^n_t} \delta_{d_j\left(\tilde{\maU}^n_t\right)}.
\end{equation}

Observe that only the order of exploration of the nodes differs in the present construction with respect to ours, this order being itself 
drawn according to a uniform choice on the set of unexplored nodes at each time in both cases. 
Therefore, as a simple consequence of the invariance in distribution of any permutation of the coordinates of $\mathbf d^n$, 
and of Proposition \ref{prop:equiv}, 

\begin{coro}
\label{cor:coincidejamming}
The sequence of measure-valued processes $\left\{\tilde\mu^n\right\}$ defined by (\ref{eq:defmutilde}) coincides in distribution with $\left\{\mu^n\right\}$, 
and so do the jamming constants of the two models.
 \end{coro}
It is significant that, although the two processes $\mu^n$ and $\tilde\mu^n$ have the same distribution, the first one is Markov but 
the second one is not, since the knowledge of the multi-graph $\tilde G$ is needed.

\subsection{Jamming constants of particular graphs}
\label{subsec:jamming}

Characterizing the jamming constant of parking problems has a long history in mathematics, see e.g. \cite{penrose}. 
One of the most studied problem in the field is the so-called random sequential absorption on discrete structures. We show hereafter how our result can 
be adapted to regular graphs, {\em i.e.} graphs with fixed degree. 
We then focus on the case of the Poisson distribution, and relate our approach to the Erd\"os-R\'enyi graph.

\paragraph{Regular graph of degree 2}
When for any $n$, the root degree distributions $\nu^n$ (and therefore, $\zeta$) is deterministic and equal to $\delta_2$, we can solve exactly the three-dimensional limiting differential system: 
\begin{equation}
\left\{\begin{array}{ll}
{d \over dt} \bar\mu_t(2) & =  - \bar\mu_t(2) (1 + 2 L(\bar\mu_t)) ;\\
& \\
{d \over dt} \bar\mu_t(1) &= - \bar\mu_t(1) (1 +  L(\bar\mu_t)) + 2 (L(\bar\mu_t)-1) \bar\mu_t(2) ;\\
&\\
{d \over dt} \bar\mu_t(0) &=  - \bar\mu_t(0)  + (L(\bar\mu_t)-1) \bar\mu_t(1) ,\label{eq:system2}
\end{array}\right.
\end{equation}
with
$$
L(\bar\mu_t)={\cro{\barmu_t,\chi^2} \over \cro{\barmu_t,\chi}}={ 4 \bar\mu_t(2) +\bar\mu_t(1) \over  2 \bar\mu_t(2) + \bar \mu_t(1) },\,t\ge 0. 
$$
After tedious but simple calculus, one obtains 
\[\left\{\begin{array}{ll}
\bar \mu_t(2) &= e^{-3t-2+2 e^{-t}};\\
 \bar \mu_t(1) &= 2(e^t-1) \bar \mu_t(2);\\
\bar \mu_t(0) &= (e^t-1)^2 \bar\mu_t(2).
\end{array}
\right.\]
Therefore the jamming constant for $\delta_2$ equals 
$$c_{\delta_2}=\int_{0}^\infty \sum_{i=0}^2 \bar \mu_t(i)dt=\int_{0}^\infty e^{-t-2+2 e^{-t}}dt= {1-e^{-2} \over 2},$$
which coincides with the famous R\'enyi parking constant of $\Z$ (see \cite{finch} Section 5.3.1). 
Informally, this can be explained by the fact that in the resulting configuration graph with distribution $\delta_2$,
only cycles of different sizes can appear while the number of small cycles (say smaller than a given constant) present in the resulting random graph is going to be very small compared
to $n$ with overwhelming probability. We do not go further into these calculations, which are beyond the scope of this paper.

\paragraph{Regular graphs for $d \ge 3$} 
More generally, jamming constants of regular graphs for $d \ge 3$ have been shown to be
 asymptotically equivalent to  ${1 - (d-1)^{-{2 \over d-2}}\over 2}$, for a random graph $\text{{\sc CM}}(n,\mathbf d^n)$ of degree distribution $\nu^n=\delta_d$ (i.e. 
$\mathbf d^n=(d,d,...,d)$), see \cite{wormaldDF}. 
 
 Let us quickly reformulate in our terminology the approach in \cite{wormaldDF}, and compare it with ours. The algorithm of {\em random pairing} as termed in \cite{wormaldDF} is similar to the one that is presented here, the main difference in the construction of $\text{{\sc CM}}(n, \mathbf d^n)$ 
being that at any time, all unexplored nodes have exactly $d$ available half-edges. This is obtained as follows: 
each time a new active node is selected (uniformly at random among all unexplored nodes), its $d$ half-edges are matched uniformly at random with $d$ other ones, exactly as we do so. At this point, according to the method of deferred decisions (see \cite{Knuth}), and unlike our construction, the connectivity of the neighbors of the new active node with the rest of the graph is not completed yet. At the following instant, a new active node is selected among all unexplored 
nodes, and the same procedure is reiterated on and on, until there is no more unexplored node. 
It is important to observe that, at the end of this algorithm, 
only the edges between active and blocked nodes have been build in the associated graph, and that the blocked nodes may have unmatched half-edges (precisely as many as $d$ minus the number of edges they share with their active neighbor). Then, the associated graph is completed by creating edges between blocked nodes in arbitrary order,  following uniform choices of half-edges (which does not change the jamming constant of the graph). Under such a dynamics, it appears clearly that the process $\left((J^n_t,U^n_t)\right)$ is Markov (recall the notation \ref{eq:defJn}). In particular, using the fact that the possible neighbors of the new active node must have $d$ available half-edges, it is easy to observe that whenever a node becomes active in $[t,t+h)$ we have that 
\begin{equation}
\label{eq:Wormald}
\esp{U^n_{t+h}-U^n_t \mid (J^n_t,U^n_t)}=-1-d{U^n_t \over n-2J^n_t},
\end{equation}
which leads to a simple one-dimensional asymptotic ODE for $t \mapsto U^n_t$, that is solved explicitly. 

However, for a general degree distribution the relation (\ref{eq:Wormald}) no longer holds, and the process $\left((J^n_t,U^n_t)\right)$ is no longer Markov. 
In fact, one can observe that the measure-valued process $\left(\mu^n_t\right)$ itself is not Markov, 
when constructing the graph as is done in \cite{wormaldDF}, since 
one needs to know which nodes are active and which are blocked in the associated graph, to complete the connectivity of the new active node at each instant - which is why we need to complete the neighboring between blocked nodes at each step. 

The price to pay for working in such generality, is that, due to a more intricate dynamics we do not obtain a closed-form formula for the function $t\mapsto \bar u_t$ in the particular case of regular graphs. However, though an exact computation of $\left(\bar u_t\right)$ becomes more and more involved when the degree $d$ grows, we can easily retrieve the asymptotic value of  
\cite{wormaldDF} by solving numerically the system corresponding to (\ref{eq:system2}), as is shown for $d=3$ and $4$ in Table \ref{table:jamming}.

\paragraph{The Poisson distribution}
\label{subsec:ER}

In the case where the asymptotic initial empirical degree distribution is Poisson, 
we obtain a closed-form expression for the function $t \mapsto \barmu_t$.

\begin{prop}
\label{pro:Poisson}
If $\zeta$ is a Poisson distribution with parameter $p$ (we denote $\zeta=\mathcal P(p)$), then
$$\barmu_t(i)= v_t {(p v_t)^i \over i!} \exp(- p v_t),\,t\ge 0,\,i\in\N,$$
where $v$ is the solution of the differential equation 
$$ {\dot v}= - v( 1 + p v).$$
Moreover the jamming constant reads 
 $$c_{\mathcal{P}(p)}=\frac{1}{p}{\log(1+p)}.$$
\end{prop}

\bp
Let us define the following Poisson measures for all $t$,   
$$\kappa_t(i) =  {(p v_t)^i \over i!} \exp(- p v_t),\,i\in\N.$$
For all $t\ge 0$, using the definition of $v_t$ and $\barmu_t$ we obtain 
\begin{align*}
{d \over dt} \cro{\barmu_t,\phi} &= \dot v_t  \exp(- p v_t) \left[   \sum_{i \ge 0} (i+1) {p^i \over i!}v_t^i \phi(i) - p   \sum_{i \ge 0}  {p^i \over i!}v_t^{i+1} \phi(i) \right],\\
 &=- \lambda v_t (1 + p v_t) \left[   \cro{\kappa_t,(1+\chi) \phi}  -  \sum_{i \ge 0}  (i+1) \kappa_t(i+1) \phi(i)\right] ,\\ 
 &=-    (1 + p v_t) \left[   \cro{\barmu_t,(1+\chi) \phi}  - \sum_{j \ge 1}  j {\barmu_t(j)} \phi(j-1)\right] ,\\ 
 &=-    \left[   \cro{\barmu_t,(1+\chi) \phi} + p v_t \cro{\barmu_t, \chi  \Delta\phi}  + p v_t  \cro{\barmu_t, \phi}   - \sum_{j \ge 1}  j \barmu_t(j) \phi(j-1) \right] ,\\ 
 &= -    \left[   \cro{\barmu_t,(1+\chi) \phi} + \left({ \cro{\barmu_t,\chi^2} \over \cro{\barmu_t,\chi}}-1\right)   \cro{\barmu_t, \chi  \Delta\phi}\right],
\end{align*}
where we use in the last identity that 
$$ p v_t = { \cro{\barmu_t,\chi^2} \over \cro{\barmu_t,\chi}}-1$$
and 
$$ p v_t   \cro{\barmu_t, \phi}=  \sum_{i \ge 0} (i+1) {p^{i+1} \over (i+1)!}v_t^{i+1} \exp(- p v_t) \phi(i)=    \sum_{j \ge 1} j {p^{j} \over j!}v_t^{j}\exp(- p v_t) \phi(j-1)= \sum_{j \ge 1}  j \barmu_t(j) \phi(j-1).$$
The jamming constant then simply follows by integration of $t \mapsto v_t$.

\ep

The asymptotic empirical degree measure of unexplored nodes thus 
turns out to be an inhomogeneous Poisson measure. A typical example where Proposition \ref{pro:Poisson} applies is that of a binomial degree distribution $\nu^n=\mbox{Bin}\,(n-1,p^n)$ for all $n$ (as is the case for Erd\"os-R\'enyi graphs), where $np^n \tend p$. 

It is then interesting to observe that we retrieve the jamming constant of the 
Erd\"os-R\'enyi graph of parameter $p$, which was given in \cite{mcdiarmid}. Notice however that we do {\em not} construct a proper 
Erd\"os-R\'enyi graph and in fact, no uniform construction based on a prescribed degree distribution can do so, since the independence assumption for the existence of the various edges cannot be fulfilled.

%
Table \ref{table:jamming} gathers in the middle column, our results for several jamming constants (exact for $\zeta=\mathcal P(1)$ and $\delta_2$, and numerically computed 
for exponential, $\delta_3$ and $\delta_4$).   
 \newcolumntype{d}{D{.}{.}{1}}
 \begin{table}[h!]
 \centering 
 \begin{tabular}{ccc} 
 \hline  
 Degree distribution & JC of Random Graphs &  JC of specific deterministic Graphs \\ [1.5ex] 
 \hline  
 Exponential $(1/2)$ &  0.7599203270 & \\
 \hline
 Poisson $(1)$ & $\log(2)=0.6931472$  &\\ 
 \hline
 $\delta_2$ & ${1-e^{-2} \over 2}=0.4323323583$  & ${1-e^{-2} \over 2}$ ($\mathbb Z$)  \\
 \hline
 $\delta_3$ & 3/8 &  0.37913944 (Honeycomb)\\
 \hline
 $\delta_4$ & 1/3 &  0.3641323  ($\mathbb Z^2$) \\
 \hline
 \end{tabular}
 \caption{Jamming constants for different degree distributions and their counterparts on deterministic graphs (simulation values for deterministic graphs are taken from \cite{rsa-sim}). }
 \label{table:jamming} 
 \end{table}

\newpage

\section{Proof of Theorem \ref{theo.main}}
\label{sec:proof}
\noindent From (\ref{eq:condinit}), there exists two real numbers $\alpha>0$ and $M>1$ such that 
\begin{align}
\cro{\zeta,\chi^6}&<M;\label{eq:defM}\\
\cro{\zeta,\chi}&>\alpha,\label{eq:defalpha}
\end{align}
where $\zeta$ is defined by (\ref{eq:convinit}). 
Let us define the following subsets of $\M_F(\N)$:
\begin{equation}
\label{eq:defMalphaM}
\M_{\alpha,M}=\biggl\{\mu\in\M_F(\N);\,\cro{\mu,\mathbf 1}\vee \cro{\mu,\chi^6} <M\mbox{ and }\cro{\mu,\chi}>\alpha\biggl\}
\end{equation}
and for all $n \in \N^*$, 
\begin{equation*}
n\M_{\alpha,M}=\biggl\{\mu\in\M_F(\N);\,{1\over n}\mu \in \M_{\alpha,M}\biggl\}.
\end{equation*}
The strategy of the proof is as follows. After showing that the generator can be approximated for large $n$ by the (non-linear) application $\Psi$ (at least on 
the subset $n\M_{\alpha,M}$),
we show uniqueness of the possible limit  - which is solution of the deterministic system of equations (\ref{eq:limflu}) - and finally the convergence in probability towards this limit.

\input{generator_approximations3.tex}

\subsection{Tightness}

\begin{prop}
\label{pro:tightness}
The sequence of measure-valued processes $\left\{\proc{\barmun}\right\}_{\N^*}$ is tight in $D([0,T],\mathcal M_F(\mathbb N))$.
\end{prop}

\bp
First fix $t\ge 0$. 
It is clear that the family of random measure $\suite{\barmun_{t}}$ is tight. 
Indeed we have almost surely for any finite subset $A$ of $\N$, $ \mu^n_t(A) \le \mu^n_0(A)$ for all $\N$. 
Hence, the family of random variables $\suite{\bar \mu^n_t(A)}$ is tight for all such $A$, which implies in turn
that the family of random measures $\suite{\barmun_t}$ is tight (see Lemma 14.15 in \cite{kallenberg}).

Therefore, from Roelly's criterion \cite{roelly}, it suffices to show that $\suite{\proc{\cro{\bar \mu_t^n,\phi}}}$ is tight in
$D\left([0,T],\R\right)$ for all $\phi \in \mathcal B_b$. For this, we exploit the semi-martingale decomposition \eqref{eq:mart} as in Joffe and M\'etivier (Corollary 2.3.3 in \cite{joffe_metivier}) and apply Rebolledo-Aldous's  criterion (\cite{aldous_tightness}) for the finite variation part
${1\over n}\int_0^.\mathscr Q\Pi_\phi\left(\mu^n_s\right)\,ds$ and the quadratic variation process $\langle\!\langle \bar M^{n}(\phi)\rangle\!\rangle$.
In detail, we aim at showing that for all $\epsilon>0$ and $\eta>0$, there exists $\delta>0$ and $n_0$ such that
\begin{align}
\sup_{n\geq n_0} &\pr{\left|{1\over n}\int_{\sigma_n}^{\tau_n}\mathscr Q\Pi_\phi\left(\mu^n_s\right)\,ds \right|\geq \eta}\leq \epsilon;
\label{eq:aldous1}\\
\sup_{n\geq n_0} &\pr{\left|\langle\!\langle \bar M^{n}(\phi)\rangle\!\rangle_{\tau_n} -\langle\!\langle \bar M^{n}(\phi)\rangle\!\rangle_{\sigma_n}\right|\geq \eta}
\leq \epsilon,\label{eq:aldous2}
\end{align}
for any two sequences $\suite{\tau_n}$ and $\suite{\sigma_n}$ of stopping times such that $\tau_n<\sigma_n<\tau_n+\delta $ for all $n\in\N^*$.
First, it readily follows from (\ref{eq:decompgenerator}) and (\ref{eq:mart}) that for all such $n\in\N^*$,
\begin{multline}
\pr{\left|{1\over n}\int_{\sigma_n}^{\tau_n}\mathscr Q\Pi_\phi\left(\mu^n_s\right)\,ds\right|\geq \eta}\\
\shoveleft{\leq \frac{\esp{\dfrac{1}{n}\displaystyle\int_{\sigma_n}^{\tau_n}\left|\mathscr Q\Pi_\phi\left(\mu^n_s\right)\right|\,ds}}{\eta}}\\
\shoveleft{\le \frac{\lambda}{\eta}\mathbf E\Biggl[\int_{\sigma_n}^{\tau_n} \cro{\barmun_s,\phi}\,ds +\int_{\sigma_n}^{\tau_n} \cro{\bar \A^n(\mu^n)_s,\phi}\,ds}\\
\shoveright{+\int_{\sigma_n}^{\tau_n}\cro{\bar \B^n(\mu^n)_s,\phi}\,ds + \int_{\sigma_n}^{\tau_n}\cro{\bar \maC^n(\mu^n)_s,\phi} ds\Biggl]}\\
\shoveleft{\le \frac{\lambda\parallel \phi \parallel}{\eta}\esp{\int_{\sigma_n}^{\tau_n}\sum_{k}\barmun_s(k)\biggl\{1+\mathbb E_Y\left[\cro{Y(\mu^n_s),\mathbf 1} \mid \mu^n_s,k\right]+\mathbb E_Y\left[\cro{Y(\mu^n_s),\chi} \mid \mu^n_s,k\right]\biggl\}\,ds}}\\
\le \frac{\lambda\parallel \phi \parallel}{\eta}\esp{(\tau_n-\sigma_n)\biggl(\cro{\barmun_0,\mathbf 1}+\cro{\barmun_0,\chi}+
\cro{\barmun_0,\chi}\cro{\barmun_0,\chi^2}\biggl)}\\
\le \frac{\lambda\parallel \phi \parallel\left(2M+M^2\right)}{\eta}\esp{\tau_n-\sigma_n},
\label{eq:tightness1}
\end{multline}
using successively Markov's inequality and the hard bounds in the proofs of Lemmas \ref{lemma:convAn} and \ref{lemma:convBn}.
All the same, using (\ref{eq:crochet}) we obtain that for some constant $C$,
\begin{multline}
\pr{\left|\langle\!\langle \bar M^{n}(\phi)\rangle\!\rangle_{\tau_n} -\langle\!\langle \bar M^{n}(\phi)\rangle\!\rangle_{\sigma_n}\right|\geq \eta}\\
\shoveleft{\le \frac{C\lambda\parallel \phi \parallel^2}{n\eta}
\mathbb E\left[\int_{\sigma_n}^{\tau_n}\sum_{k}\barmun_s(k)
\Biggl\{1+\left(\mathbb E_Y\left[\cro{Y(\mu^n_s),\mathbf 1} \mid \mu^n_s,k\right]\right)^2\right.}\\
\shoveright{\left.+\left(\mathbb E_Y\left[\cro{Y(\mu^n_s),\chi} \mid \mu^n_s,k\right]\right)^2\Biggl\}\,ds\right]}\\
\le \frac{C\lambda\parallel \phi \parallel^2}{n\eta}
\esp{(\tau_n-\sigma_n)\biggl(\cro{\barmun_0,\mathbf 1}+\cro{\barmun_0,\chi^2}+
\cro{\barmun_0,\chi^2}\left(\cro{\barmun_0,\chi^2}\right)^2\biggl)}\\
\le \frac{C\lambda\parallel \phi \parallel^2\left(2M+M^3\right)}{n\eta}\esp{\tau_n-\sigma_n}.
\label{eq:tightness2}
\end{multline}

Clearly, from (\ref{eq:tightness1}) and (\ref{eq:tightness2}) we can choose $\delta$ small enough so that (\ref{eq:aldous1}) and (\ref{eq:aldous2}) hold for
$n_0=1$, which concludes the proof.

\ep

\subsection{Uniqueness}

Let $\proc{\bar \mu}$ and $\proc{\bar \nu}$, two solutions of (\ref{eq:limflu}) in $\mathcal C\left(\R+,\M_F(\N)\right)$, with the same initial condition 
$\zeta$, and let us denote
 $\gamma_t=\bar \mu_t-\bar \nu_t$ for all $t$. Denote also for all $\beta$ such that $0<\beta<\alpha$,
\begin{align*}
t_\beta^{\mu}&=\sup\Bigl\{t\ge 0;\,\cro{\barmu_t,\chi}>\beta\Bigl\};\\
t_\beta^{\nu}&=\sup\Bigl\{t\ge 0;\,\cro{\bar\nu_t,\chi}>\beta\Bigl\}.
\end{align*}
Recall, that in view of (\ref{eq:defM}) we have that
$$  \cro{\barmu_t,\chi^5}\vee\cro{\bar\nu_t,\chi^5}\le M,\,t\ge 0.$$
Thus, denoting again for all $\mu \in \M_F(\N)$, $L(\mu)={\cro{\mu,\chi^2} \over \cro{\mu,\chi}}$ we have that
$$L\left(\barmu_t\right)\vee L\left(\bar\nu_t\right) \le {M\over \beta},\,t\in \left[0,t^\mu_{\beta}\wedge t^\nu_{\beta}\right].$$

\bigskip
\noindent
Let for all $t\ge 0$,
$$\Gamma_t=\sum_{i \ge 0} i^8 \gamma_t(i)^2.$$
Fix $\beta$ such that $0<\beta<\alpha$ and $t\in \left[0,t^\mu_{\beta}\wedge t^\nu_{\beta}\right]$. From (\ref{eq:limflu}), we have that
\begin{align}
{d \over dt} \Gamma_t&=\sum_{i \ge 0}i^8 {d \over dt}\gamma_t(i)^2\nonumber\\
&= 2\sum_{i \ge 0} i^8 \gamma_t(i){d \over dt}\gamma_t(i)\nonumber\\
&= 2\lambda \sum_{i \ge 0} \biggl\{(i+1)\left(L(\bar \mu_t)-1\right)\gamma_t(i)\gamma_t(i+1) -i^8\left(1+ i L(\bar \mu_t)\right)\left(\gamma_t(i)\right)^2\biggl\}\nonumber\\
&\phantom{2\lambda \sum_{i \ge 0} \biggl\{(i+1)}-2\lambda \sum_{i \ge 0}\Bigl\{\left(L(\bar \nu_t)-L(\bar \mu_t)\right)  i^8 \gamma_t(i)(i \bar \nu_t(i)-(i+1)\bar \nu_t(i+1))\biggl\}\nonumber\\
&= A_t + B_t.\label{eq:in0}
\end{align}
We first deal with the term $A_t$. Using that
\begin{equation*}
(i+1)\gamma_t(i+1)\gamma_t(i) = {1\over 2}  (i+1)\gamma_t(i+1)^2 + {1\over 2} (i+1)\gamma_t(i)^2
- {1\over 2}\Bigl(\sqrt{i+1}\gamma_t(i+1) - \sqrt{i+1}\gamma_t(i)\Bigl)^2,
\end{equation*}
we get
 \begin{align*}
 \sum_{i\ge 0} i^8(i+1)\gamma_t(i+1)\gamma_t(i) &\le  {1\over 2}\sum_{i \ge 0} i^8(i+1) \gamma_t(i+1)^2  + {1\over 2}\sum_{i \ge 0} i^8(i+1)\gamma_t(i)^2 \\
 \nonumber
 &\le {1\over 2}\sum_{i \ge 0} (i+1)^9 \gamma_t(i+1)^2  + {1\over 2}\sum_{i \ge 0} i^9 \gamma_t(i)^2 + {1\over 2}\sum_{i \ge 0} i^8 \gamma_t(i)^2\\
 &\le \sum_{i \ge 0} i^9 \gamma_t(i)^2 + {1\over 2}\sum_{i \ge 0} i^8 \gamma_t(i)^2.
\end{align*}
Hence, as $L(\barmu_t) \ge 1$,
\begin{align}
A_t&=2\lambda\sum_i  \biggl\{i^8 \left(L(\bar\mu_t)-1\right)(i+1)\gamma_t(i+1)\gamma_t(i)+ i^8\left(1+iL\left(\bar\mu_t\right)\right)\left(\gamma_t(i)\right)^2\biggl\}\nonumber\\
&\le 2\lambda(L(\bar\mu_t) -1)\left( \sum_{i \ge 0}   i^9 \gamma_t(i)^2  + {1\over 2} \sum_{i \ge 0}  i^8 \gamma_t(i)^2 \right) - 2\lambda L(\bar\mu_t)  \sum_{i \ge 0}  i^9  \gamma_t(i)^2 \nonumber\\
&\le  \lambda L\left(\barmu_t\right)\sum_i i^8  \gamma_t(i)^2\nonumber\\
&\le {\lambda M\over \beta}\Gamma_t.\label{eq:in2}
\end{align}

\bigskip
\noindent Let us now deal with $B_t$. First observe that
\begin{multline*}
|L(\bar \mu_t)- L(\bar \nu_t)|
\\
\shoveleft{
=\left|{ \sum_i i^2\bar \mu_t(i) \sum_i i\bar \nu_t(i) -\sum_i i\bar \mu_t(i)\sum_i i^2\bar \nu_t(i) \over \sum_i i\bar \mu_t(i)\sum_i i\bar \nu_t(i)} \right|}\\
= \left|{ \sum_i i^2\gamma_t(i) \sum_i i\bar \nu_t(i) +\sum_i i^2\bar \nu_t(i)\sum_i i\bar \nu_t(i) -\sum_i i\gamma_t(i)\sum_i i^2\bar \nu_t(i) - \sum_i i^2\bar \nu_t(i)\sum_i i\bar \nu_t(i) \over \sum_i i\bar \mu_t(i)\sum_i i\bar \nu_t(i)} \right| \\
\le { \sum_i i^2|\gamma_t(i)| \sum_i i\bar \nu_t(i) +\sum_i i |\gamma_t(i)|\sum_i i^2\bar \nu_t(i) \over \sum_i i\bar \mu_t(i)\sum_i i\bar \nu_t(i)}
\le   {2M \over \beta^2} \sum_i i^2|\gamma_t(i)|.
\end{multline*}
Hence,
\begin{align}
\left|L(\bar \mu_t)- L(\bar \nu_t)\right| &\le {2M \over \beta^2}\sum_i i^2 |\gamma_t(i)|,\nonumber\\
 &\le {2M \over \beta^2}\left(\sup_i i^4 |\gamma_t(i)|\right) \sum_i i^{-2}\nonumber\\
 &\le {M\pi^2 \over 3\beta^2}\left(\sup_i i^8 |\gamma_t(i)|^2\right)^{1/2}\nonumber\\
 &\le {M\pi^2 \over 3\beta^2}\left(\sum i^8 |\gamma_t(i)|^2\right)^{1/2}.\label{in4}
\end{align}
On the other hand, using Cauchy Schwartz inequality,
\begin{align}
\sum_i  \bar \nu(i)i^9 |\gamma_t(i)| &=  \sum_i i^5 \bar \nu_t(i) i^4 |\gamma_t(i)|\nonumber\\
&\le \left( \sum_i i^{10} \bar \nu_t(i)^2\right)^{1/2}  \left(\sum_i i^8 \gamma_t(i)^2\right)^{1/2}\nonumber\\
&\le  \left(\sum_i i^{5} \bar \nu_t(i)\right) \left(\sum_i i^8 \gamma_t(i)^2\right)^{1/2}\nonumber\\
&\le M  \left(\sum_i i^8 \gamma_t(i)^2\right)^{1/2},
\label{in5}
\end{align}
where in the second inequality we used the fact that for a finite series $x$ having positive terms,
$$\sum_{i\in\N} x_i^2 \le \left(\sum_{i\in\N} x_i\right)^2.$$
All the same, we have that
\begin{align}
\sum_i  \bar \nu_t(i+1)(i+1)i^8 |\gamma_t(i)| &=  \sum_i (i+1)^5 \bar \nu_t(i+1) i^4 |\gamma_t(i)|\nonumber\\
&\le \left( \sum_i (i+1)^{10} \bar \nu_t(i+1)^2\right)^{1/2}  \left(\sum_i i^8 \gamma_t(i)^2\right)^{1/2}\nonumber\\
&\le \left( \sum_i i^{10} \bar \nu_t(i)^2\right)^{1/2}  \left(\sum_i i^8 \gamma_t(i)^2\right)^{1/2}\nonumber\\
&\le M  \left(\sum_i i^8 \gamma_t(i)^2\right)^{1/2}.
\label{in5bis}
\end{align}
Hence, using (\ref{in4}), (\ref{in5}) and (\ref{in5bis}) we obtain that
\begin{align}
B &\le  2\lambda |L(\bar \mu_t)- L(\bar \nu_t)|\sum_i  \bar \nu_t(i)i^9 |\gamma_t(i)|  + \bar \nu_t(i+1)i^8(i+1) |\gamma_t(i)| \nonumber\\
&\le   {4\lambda M^2\pi^2 \over 3\beta^2}\left(\sum_i i^8 |\gamma_t(i)|^2\right).\label{in6}
\end{align}

Finally, using (\ref{eq:in2}) and (\ref{in6}) in $(\ref{eq:in0})$, we obtain that for some positive constant $C$, for all $t\le t^{\mu}_\beta\wedge t^{\nu}_\beta$,
\begin{eqnarray*}
{d \over dt} \Gamma_t \le C \Gamma_t.
\end{eqnarray*}
Since $G(0)=0$, $G$ is a positive function and $t^{\mu}_\beta\wedge t^{\nu}_\beta>0$, this shows using Gronwall's Lemma that $\Gamma_t=0$ for all such $t$. Therefore, $t^{\mu}_\beta=t^{\nu}_\beta=:t_\beta$, and $\barmu_t$ and $\bar\nu_t$ coincide up to $t_\beta$. In other words there is at most one
solution to (\ref{eq:limflu}) up to time $t_\beta$. Since this is true for all $\beta$, and since the only solution $\barmu$ is such that
$t\to \cro{\barmu_t,\chi}$ is continuous, there is at most one solution up to the (positive, in view of (\ref{eq:condinit})) instant
\begin{equation}
\label{eq:deft0}
t_0=\sup\left\{t_\beta;\,\beta>0\right\}.
\end{equation}
The proof of uniqueness is completed by noticing that whenever $t_0<\infty$, the only solution $\barmu$ to (\ref{eq:limflu}) can be extended uniquely after $t_0$, as follows:
\begin{equation}
\label{eq:aftert0}
\barmu_t=\mu_{t_0}(0)e^{-\lambda(t-t_0)}\delta_0,\,t\ge t_0.
\end{equation}


\subsection{Convergence}
\label{subsec:cv}
Recall that $\zeta\in \M_f(\N)$, $M>1$ and $\alpha>0$ are respectively defined by (\ref{eq:convinit}), (\ref{eq:defM}) and (\ref{eq:defalpha}). 
Fix $T>0$.   

\paragraph*{$\underline{\mbox{Convergence before reaching a given positive threshold}}$}
\text{ }\\

\noindent We first prove the following result.
\begin{prop}
\label{prop:horizonbefore}
For all $\phi \in \mathcal B_b$, 
\begin{equation*}
\sup_{t\in[0,\t2a]}\left|\cro{\barmun_t,\phi} -\cro{\barmu_t,\phi}\right| \overset{(\mathcal P)}{\tend} 0,
\end{equation*}
where $\proc{\bar\mu}$ is the unique solution of (\ref{eq:limflu}) on $\C\left([0,\t2a],\M_F(\N)\right)$. 
\end{prop}
\bp
For all $n\in\N^*$, we define the stopping times
$$\tau^{n}_\alpha =\sup\left\{t \ge 0;\,  \cro{\mu^n_t,\chi}\ge n\alpha \right\}$$
and denote for all $t$,
$$\barmunt_t=\barmun_{t\wedge \tau^n_\alpha}.$$
Clearly, from Proposition \ref{pro:tightness} the sequence of stopped processes $\left\{\barmunt\right\}$ is also relatively compact for the topology of weak convergence, and 
we let $\barmus$ be a sub-sequential limit. Define
\begin{equation*}
\T2a=\sup\Bigl\{t \ge 0;\,\cro{\barmus_t,\chi}>2\alpha\Bigl\}.
\end{equation*}
Notice that the process $\barmus$ and hence the instant $\T2a$, are {\em a priori} random. 
Fix $\phi\in\mathcal B_b$ throughout the proof. Let $n\in\N^*$. 
We have for all $t\in [0,T]$,
\begin{multline}
\int_0^{t\wedge\T2a\wedge\tau_\alpha^n}\cro{\Psi\left(\barmunt\right)_s,\phi}\,ds
=\left(\int_0^{t\wedge\T2a}\cro{\Psi\left(\barmunt\right)_s,\phi}\,ds\right)\ind_{\{\tau_\alpha^n > \T2a\wedge T\}}\\
  +\left(\int_0^{t\wedge\tau_\alpha^n}\cro{\Psi\left(\barmunt\right)_s,\phi}\,ds\right)\ind_{\{\tau_\alpha^n \le \T2a\wedge T\}}\label{eq:CMT}.
\end{multline}
On the one hand, in view of Lemma A.5 of (\cite{DDMT12}), the following map is continuous for the Skorokhod topology:
\[\left\{\begin{array}{ccc}
 \D\left([0,T],\M_{\alpha,M}\right) &\longrightarrow &\D\left([0,T],\R\right)\\
 \mu &\longmapsto &\cro{\mu_.,\chi},
             \end{array}
             \right.\]
and so does
\[\left\{\begin{array}{ccc}
 \D\left([0,T],\R\right) &\longrightarrow &\R\\
 x_. &\longmapsto &\inf_{t\in[0,T]}x_t.
             \end{array}
             \right.\]
Therefore, from the Continuous Mapping Theorem (\cite{Bil68}), along the latter subsequence the following convergence in distribution holds:
$$\inf_{t\in[0,T]} \cro{\barmunt_{t\wedge \T2a},\chi} \Rightarrow \inf_{t\in[0,T]} \cro{\barmus_{t\wedge \T2a},\chi}.$$
Hence, from Fatou's Lemma,
\begin{equation}
1=\pr{\inf_{t\in [0,T]} \cro{\barmus_{t\wedge \T2a},\chi} > \alpha} \le \lim_{n\rightarrow \infty}
\pr{\inf_{t\in [0,T]} \cro{\barmunt_{t\wedge \T2a},\chi} > \alpha}
\le \lim_{n\rightarrow \infty}\pr{T \wedge \T2a <\tau^n_\alpha}.
\label{eq:Fatou}
\end{equation}
\noindent Therefore,
\begin{equation}
\left(\int_0^{t\wedge\tau_\alpha^n}\cro{\Psi\left(\barmunt\right)_s,\phi}\,ds\right)\ind_{\{\tau_\alpha^n \le \T2a\wedge T\}}
\tendP 0.
\label{eq:rhs}
\end{equation}
\medskip

\noindent
Now, it follows again from Lemma A.5 of (\cite{DDMT12}) that the following mappings are continuous for the Skorokhod topology:
\[\left\{\begin{array}{ccc}
 \D\left([0,T],\M_{\alpha,M}\right) &\longrightarrow &\D\left([0,T],\R\right)\\
 \mu &\longmapsto &\cro{\mu_.,\chi};\\
         &    &\cro{\mu_.,\chi^2};\\
           &  &\cro{\mu_.,\left(\mathbf 1+\chi\right)\phi};\\
            & &\cro{\mu_.,\chi\Delta\phi},
             \end{array}
             \right.\]
and it is a classical result that the followings map is also continuous: 
\[\left\{\begin{array}{ccc}
 \D\left([0,T],\R\times\R^*\right) &\longrightarrow &\D\left([0,T],\R\right)\\
 (x_.,y_.) &\longmapsto &{x_. \over y_.}.\\
             \end{array}
             \right.\]
So from the Continuous Mapping Theorem, the map
\[\left\{\begin{array}{ccc}
 \D\left([0,T],\M_{\alpha,M}\right) &\longrightarrow &\D\left([0,T],\R\right)\\
 \mu &\longmapsto &\cro{\Psi(\mu)_t,\phi}\\
             \end{array}
             \right.\]
is itself continuous, and it follows from the continuity of the map
\[\left\{\begin{array}{ccc}
 \D\left([0,T],\R\right) &\longrightarrow &\C\left([0,T],\R\right)\\
 x_. &\longmapsto &\int_0^. x_s\,ds,\\
             \end{array}
             \right.\]
together with (\ref{eq:CMT}), (\ref{eq:Fatou}) and (\ref{eq:rhs}), that along the same sub-sequence
\begin{equation}
\label{eq:CMT1}
\int_0^{.\wedge\T2a\wedge \tau^n_\alpha}\cro{\Psi\left(\barmunt\right)_s,\phi}\,ds \Longrightarrow \int_0^{.\wedge\T2a}\cro{\Psi\left(\barmus\right)_s,\phi}\,ds\,\, \mbox{ in }\C\left([0,T],\R\right).
\end{equation}

\bigskip
\noindent
On the other hand, whenever $\barmun_0 \in \M_{\alpha,M}$ we clearly have that $\barmunt_{t\wedge \T2a} \in \M_{\alpha,M}\mbox{ for all }t$.
Therefore, as a consequence of Lemmas \ref{lemma:convCn}, \ref{lemma:convAn} and \ref{lemma:convBn} together with (\ref{eq:condinit}), (\ref{eq:mart}), (\ref{eq:decompPsi}) and (\ref{eq:decompgenerator}), we have a.s. for all $n\in\N^*$ and $t\ge 0$, 
\begin{multline}
\cro{\barmunt_{t\wedge \T2a},\phi}=\cro{\barmun_0,\phi}
+\left(\int_0^{t\wedge\tau^n_\alpha\wedge \T2a}\cro{\Psi\left(\barmunt\right)_s,\phi}\,ds +o^{n,\alpha}_t\right)\ind_{\{\barmun_0 \in \M_{\alpha,\M}\}}\\
+\left({1 \over n}\int_0^{t\wedge\tau^n_\alpha\wedge \T2a}\mathscr Q^n\Pi_\phi\left(\mu^n_s\right) \,ds \right)\ind_{\{\barmun_0\not\in \M_{\alpha,M}\}}+\bar M^n(\phi)_{t\wedge \tau^n_\alpha\wedge \T2a},
\label{eq:prelimit}
\end{multline}
where $o^{n,\alpha}$ is a process converging to 0 in probability and uniformly over compact sets.\\

\noindent
Now, applying Doob's inequality to the stopped martingale $\bar M^n(\phi)_{.\wedge \tau^n_\alpha\wedge \T2a}$ and using
(\ref{eq:crochet}) as in (\ref{eq:tightness2}) yields that
\begin{equation}
\label{eq:prelimit1}
\sup_{t\in[0,T]}\left|\bar M^n(\phi)_{t\wedge \tau^n_\alpha}\right| \overset{(\mathcal P)}{\tend} 0.
\end{equation}
Moreover, for all $n$, $t$ and all $\epsilon$ we have that
\begin{multline}
\pr{\sup_{t\in[0,T]}\left({1 \over n}\int_0^{t\wedge\tau^n_\alpha\wedge \T2a}\mathscr Q^n\Pi_\phi\left(\mu^n_s\right) \,ds \right)\ind_{\{\barmun_0\not\in \M_{\alpha,M}\}}>\epsilon}\\
\shoveleft{= \pr{\sup_{t\in[0,T]} {1 \over n}\int_0^{t\wedge\tau^n_\alpha\wedge \T2a}\mathscr Q^n\Pi_\phi\left(\mu^n_s\right) \,ds >\epsilon, \barmun_0\not\in \M_{\alpha,M}}}\\
\shoveleft{\le \pr{\barmun_0\not\in \M_{\alpha,M}}}\\
\shoveleft{\le \pr{\cro{\barmun_0,\chi}\le \alpha}+\pr{\cro{\barmun_0,\chi^6} \ge M}}\\
\shoveleft{\le \pr{\cro{\barmun_0,\chi}< \cro{\barmu_0,\chi}-{\cro{\barmu_0,\chi}-\alpha \over 2}}
+\pr{\cro{\barmun_0,\chi^6} > \cro{\barmu_0,\chi^6} +{M-\cro{\barmu_0,\chi^6} \over 2}}}
\tend 0,\label{eq:prelimit2}
\end{multline}
in view of (\ref{eq:convinit}) and (\ref{eq:condinit}).
Plugging (\ref{eq:convinit}) together with (\ref{eq:CMT1}), (\ref{eq:prelimit1}) and (\ref{eq:prelimit2}) into (\ref{eq:prelimit}), and using Skorokhod Representation Theorem implies that on some probability space, almost surely
\begin{equation*}
\cro{\barmus_{t\wedge\T2a},\phi}=\cro{\barmu_0,\phi}+\int_0^{t\wedge\T2a}\cro{\Psi\left(\barmus\right)_s,\phi}\,ds,\,t\ge 0.
\end{equation*}
In other words, $\barmus$ is a $\C\left([0,T],\M_F(\N)\right)$-valued process having initial deterministic value $\zeta$, 
and solving (\ref{eq:limflu}) on $\left[0,T\wedge \T2a\right]$. 
As the solution of the latter, if any, is unique, we conclude (i) that $\T2a$ is deterministic, hence (ii) that there exists a solution $\barmu$ to (\ref{eq:limflu}) on $\left[0,T\wedge\T2a\right]$, with which $\barmus$ coincides almost surely and (iii) that $\T2a=\t2a$, where 
\begin{equation}
\label{eq:deft2a}
\t2a=\sup\Bigl\{t\ge 0;\,\cro{\barmu_t,\chi}>2\alpha\Bigl\}.
\end{equation}
In particular, as the tightness of $\{\barmunt\}$ clearly implies that of $\left\{\int_0^{.\wedge\T2a\wedge \tau^n_\alpha}\cro{\Psi\left(\barmunt\right)_s,\phi}\,ds\right\}$, we deduce from (\ref{eq:CMT1}) that
\begin{equation*}
\int_0^{.\wedge\t2a\wedge \tau^n_\alpha}\cro{\Psi\left(\barmunt\right)_s,\phi}\,ds \Longrightarrow \int_0^{.\wedge\t2a}\cro{\Psi\left(\barmu\right)_s,\phi}\,ds\,\, \mbox{ in }\C\left([0,T],\R\right).
\end{equation*}
Using once again the Representation Theorem we obtain that on some probability space,
\begin{equation*}
\int_0^{.\wedge\t2a\wedge \tau^n_\alpha}\cro{\Psi\left(\barmunt\right)_s,\phi}\,ds \tend \int_0^{.\wedge\t2a}\cro{\Psi\left(\barmu\right)_s,\phi}\,ds\,\, \mbox{ a.s. in }\C\left([0,T],\R\right)
\end{equation*}
which, as the Skorokhod topology and the uniform topology coincide on $\C\left([0,T],\R\right)$ (see \cite{Bil68}, p.112), implies that
\begin{equation*}
\sup_{t\in[0,T]}\left|\int_0^{t\wedge\t2a\wedge \tau^n_\alpha}\cro{\Psi\left(\barmunt\right)_s,\phi}\,ds -\int_0^{t\wedge\t2a}\cro{\Psi\left(\barmu\right)_s,\phi}\,ds\right| \tend 0\,\, \mbox{ a.s.}
\end{equation*}
and therefore,
\begin{equation*}
\sup_{t\in[0,T]}\left|\int_0^{t\wedge\t2a\wedge \tau^n_\alpha}\cro{\Psi\left(\barmunt\right)_s,\phi}\,ds -\int_0^{t\wedge\t2a}\cro{\Psi\left(\barmu\right)_s,\phi}\,ds\right| \overset{(\mathcal P)}{\tend} 0.
\end{equation*}
As the latter holds true for all $T>0$, we obtain that
\begin{equation}
\label{eq:CMT2}
\sup_{t\in[0,\t2a]}\left|\int_0^{t\wedge \tau^n_\alpha}\cro{\Psi\left(\barmunt\right)_s,\phi}\,ds -\int_0^{t}\cro{\Psi\left(\barmu\right)_s,\phi}\,ds\right| \overset{(\mathcal P)}{\tend} 0.
\end{equation}
Consequently, from (\ref{eq:prelimit}) we obtain that for all $\epsilon>0$,
\begin{multline*}
\pr{\sup_{t\in\left[0,\t2a\right]}\left|\cro{\barmun_t,\phi}-\cro{\barmu_t,\phi}\right|>\epsilon}\\
\le \pr{\left\{\sup_{t\in\left[0,\t2a\right]}\left|\int_0^{t\wedge \tau^n_\alpha}\cro{\Psi\left(\barmunt\right)_s,\phi}\,ds -\int_0^{t}\cro{\Psi\left(\barmu\right)_s,\phi}\,ds\right|>{\epsilon\over 3}\right\}\bigcap\left\{\tau^n_\alpha> \t2a\right\}}\\
+\pr{\sup_{t\in\left[0,\t2a\right]}\left|\bar M^n(\phi)_{t\wedge \tau^n_\alpha}+o^{n,\alpha}_t\right|>{\epsilon\over 3}}
+\pr{\left|\cro{\barmunt_0,\phi}-\cro{\barmu_0,\phi}\right|>{\epsilon \over 3}}
+\pr{\tau^n_\alpha\le \t2a}.
\end{multline*}
The first term on the r.h.s. vanishes for large $n$ thanks to (\ref{eq:CMT2}), the second one from Doob's inequality, the third one in view of (\ref{eq:convinit}) and the last one from
(\ref{eq:Fatou}). This concludes the proof.
\ep

\paragraph*{$\underline{\mbox{Existence of the solution on }\R+}$}

A consequence of Proposition \ref{prop:horizonbefore} is the existence of a solution $\proc{\barmu}$ of (\ref{eq:limflu}) 
until $\t2a$. As a matter of fact, it appears clearly that the latter result can be extended to any 
$0<\beta<\alpha$. Therefore, by its continuity the solution $\proc{\barmu}$ can be extended at least until $t_0$, the 
hitting time of 0 defined by (\ref{eq:deft0}). 
The existence of the solution $\proc{\barmu}$ after $t_0$ then follows from the explicit form (\ref{eq:aftert0}).  

\paragraph*{$\underline{\mbox{Asymptotics of the mass at the origin }}$}

We now focus on the mass concentrated at 0 in the hydrodynamic limit. 
Let us first give the following result,  

\begin{lemma}
\label{lemma:masse0}
For all $t \ge \t2a$,
$$\lambda\int_{\t2a}^{t\wedge t_0} \barmu_s(1)\left({\cro{\barmu_s,\chi^2} \over \cro{\barmu_s,\chi}}-1\right) ds\le \cro{\barmu_{\t2a},\chi},$$
where $\proc{\barmu}$ is the only solution to (\ref{eq:limflu}) on $\R+$ and $\t2a$ is defined by (\ref{eq:deft2a}). 
\end{lemma}

\bp
Let $t\ge \t2a$. Plainly, applying (\ref{eq:limflu}) to $\phi:=\chi$ leads to
\begin{multline*}
\cro{\barmu_t,\chi}=\cro{\barmu_{\t2a},\chi}-\lambda\left\{\int_{\t2a}^t \cro{\barmu_s,\chi}\,ds+
\int_{\t2a}^{t\wedge t_0} \cro{\barmu_s,\chi^2}\,ds\right.\\
\left.+\int_{\t2a}^{t\wedge t_0} \cro{\barmu_s,\chi\Delta\chi}\left({\cro{\barmu_s,\chi^2} \over \cro{\barmu_s,\chi}}-1\right)\,ds\right\}.
\end{multline*}
Therefore, as $\cro{\barmu_t,\chi}\ge 0$ we obtain that
\begin{equation*}
\cro{\barmu_{\t2a},\chi}
                    \ge \lambda\int_{\t2a}^{t\wedge t_0} \cro{\barmu_s,\chi\Delta\chi}\left({\cro{\barmu_s,\chi^2} \over \cro{\barmu_s,\chi}}-1\right)\,ds.
\end{equation*}
The proof is completed by noticing that for all $s$,
$$\cro{\barmu_s,\chi\Delta\chi}=\sum_{i\in\N}\barmu_s(i)i\left(i-(i-1)\right)=\sum_{i\in\N}\barmu_s(i)i\ge \barmu_s(1).$$
\ep

\noindent Now, let for all $n\in\N^*$ and all $t\ge \t2a$,
\begin{align*}
X^{n,\alpha}_t&=\mbox{Card}\,\Bigl\{\mbox{$\sU$-vertices of degree 0 at $\t2a$ that have  become $\sA$ before time $t$}\Bigl\};\\
Y^{n,\alpha}_t&=\mbox{Card}\,\Bigl\{\mbox{$\sU$-vertices of degree $\ge 1$ at $\t2a$,}\\
              &\,\,\,\,\,\,\,\,\,\,\,\,\,\,\,\,\,\,\,\,\,\,\,\,\,\,\, \mbox{ having become of degree 0 at time $t$ and being still $\sU$ at $t$}\Bigl\}.\\
\end{align*}
Denote also $\bar X^{n,\alpha}_t={1\over n}X^{n,\alpha}_t$ and $\bar Y^{n,\alpha}_t={1\over n}Y^{n,\alpha}_t$ for all $n$ and $t$.
Then, clearly
\begin{equation}
\label{eq:evolmass0norm}
\barmun_t(0)=\barmun_{\t2a}\left(0\right)-\bar X^{n,\alpha}_t+\bar Y^{n,\alpha}_t,\,\,t \ge \t2a.
\end{equation}
As the $\sU$-vertices of degree 0 are independent of the rest of the graph and eventually all become $\sA$-vertices, it is clear that for any $n$, the process $X^{n,\alpha}$ is Markov on $\N$. It has rcll paths, its generator clearly reads
$$\tilde{\mathscr Q}F(x)=\lambda \left(\mu^n_{\t2a}(0)-x\right)$$
for all functions $F:\R \to \R$ and all $x\in\N$, so it is routine to check that for all $n\in\N^*$, for some square integrable martingale 
$\bar M^{n,\alpha}$, for all $t\ge \t2a$,
\begin{equation}
\label{eq:martX}
\bar X^{n,\alpha}_t=\lambda\int_{\t2a}^t\left(\barmun_{\t2a}(0)-\bar X^{n,\alpha}_s\right)\,ds+\bar M^{n,\alpha}_t.
\end{equation}
Therefore, with (\ref{eq:evolmass0norm}), we obtain that
\begin{align}
\barmun_t(0)&=\barmun_{\t2a}(0)-\lambda\int_{\t2a}^t \left(\barmun_{\t2a}(0)-\bar X^{n,\alpha}_s\right)\,ds-\bar M^{n,\alpha}_t+\bar Y^{n,\alpha}_t\nonumber\\
            &=\barmun_{\t2a}(0)-\lambda\int_{\t2a}^t \left(\barmun_{s}(0)-\bar Y^{n,\alpha}_s\right)\,ds-\bar M^{n,\alpha}_t+\bar Y^{n,\alpha}_t.
\label{eq:afterhorizon0}
\end{align}
It also readily follows once again by Doob's inequality that the martingale term vanishes uniformly in $L^2$, and in particular, that for all $t\ge \t2a$,
\begin{equation}
\label{eq:convM}
\sup_{s\in [\t2a,t]}  \left|\bar M^{n,\alpha}_s\right| \overset{(\mathcal P)}{\tend} 0.
\end{equation}

\noindent
On the other hand, applying (\ref{eq:limflu}) to $\phi\equiv \ind_0$, and observing that for all $i\in\N$,
$$\chi(i)\Delta\ind_0(i)=-\ind_1(i),$$
we get that
$$\barmu_t(0)=\barmu_{\t2a}(0)-\lambda\left\{\int_{\t2a}^t \barmu_s(0)\,ds-\int_{\t2a}^{t\wedge t_0} \barmu_s(1)\left({\cro{\barmu_s,\chi^2} \over \cro{\barmu_s,\chi}}-1\right)\,du\right\},\,t\ge \t2a.$$
Combining this with (\ref{eq:afterhorizon0}) leads for all $n$ to
\begin{multline}
\left|\barmun_t(0)-\barmu_t(0)\right|\le \left|\barmun_{\t2a}(0)-\barmu_{\t2a}(0)\right|+\lambda \int_{\t2a}^t \left|\barmun_s(0)-\barmu_s(0)\right|\,ds
+\left|\bar M^{n,\alpha}_t\right|\\
+\bar Y^{n,\alpha}_t+\lambda\int_{\t2a}^t \bar Y^{n,\alpha}_s\,ds
+\lambda\int_{\t2a}^t \barmu_s(1)\left({\cro{\barmu_s,\chi^2} \over \cro{\barmu_s,\chi}}-1\right)\,ds.\label{eq:afterhorizon1}
\end{multline}
But by its very definition, we have that for all $s \ge \t2a$ and all $n$,
$$\bar Y^{n,\alpha}_s \le \barmun_{\t2a}(1)\le\cro{\barmun_{\t2a},\chi}.$$
Plugging this together with Lemma \ref{lemma:masse0} in (\ref{eq:afterhorizon1}) yields to
\begin{multline*}
\left|\barmun_t(0)-\barmu_t(0)\right|\le \left|\barmun_{\t2a}(0)-\barmu_{\t2a}(0)\right|+\lambda \int_{\t2a}^t \left|\barmun_s(0)-\barmu_s(0)\right|\,ds
+\left|\bar M^{n,\alpha}_t\right|\\
+\left(1+\lambda (t-\t2a)\right)\cro{\barmun_{\t2a},\chi}
+\cro{\barmu_{\t2a},\chi}.
\end{multline*}
Therefore, from Gronwall's Lemma we conclude that for all $n\in \N^*$ and all $t \ge \t2a$, 
\begin{equation}
\label{eq:afterhorizon2}
\left|\barmun_t(0)-\barmu_t(0)\right|
\le \Bigl(\left|\barmun_{\t2a}(0)-\barmu_{\t2a}(0)\right|+\left|\bar M^{n,\alpha}_t\right|
+\left(1+\lambda (t-\t2a)\right)\cro{\barmun_{\t2a},\chi}+2\alpha \Bigl)e^{\lambda (t-\t2a)}.
\end{equation}

\paragraph*{$\underline{\mbox{Proof of Theorem \ref{theo.main}}}$}
We are now in position to prove Theorem \ref{theo.main}. 
Fix $T>0$. First consider a function $\phi \in \maB_b$ such that $\phi(0) \ne 0$. Let $\epsilon>0$.
We can chose $\alpha$ small enough in (\ref{eq:defalpha}), and small enough positive numbers $\delta$, $\eta$ and $\xi$ so that 
\begin{equation}
\Bigl(\delta+\eta+\left(1+\lambda (T-\t2a)^+\right)\left(\alpha+\xi\right)+2\alpha \Bigl)e^{\lambda (T-\t2a)^+} <{\varepsilon \over 2\mid \phi(0) \mid};
\label{eq:condalpha1}
\end{equation}
\begin{equation}
4\alpha+\xi <{\varepsilon \over 2\parallel \phi\parallel}.\label{eq:condalpha2}
\end{equation}
First, if $T\le \t2a$, Proposition \ref{prop:horizonbefore} trivially implies that
\begin{equation*}
\pr{\sup_{t\in[0,T]}\left|\cro{\barmun_t,\phi} -\cro{\barmu_t,\phi}\right| >\epsilon}\tend 0.
\end{equation*}
If $T>\t2a$, define the following events  for all $n\in\N^*$, 
\begin{align*}
\Omega^n_{\alpha,\delta}&=\Bigl\{\sup_{t\in (0,\t2a]}\left|\barmun_{t}(0)-\barmu_{t}(0)\right|\le \delta\Bigl\};\\
\Omega^n_{\alpha,\eta}&=\Bigl\{\sup_{t\in (\t2a,T]}\left|\bar M^{n,\alpha}_t\right| \le \eta\Bigl\};\\
\Omega^n_{\alpha,\xi}&=\Bigl\{\sup_{t\in (0,\t2a]}\left|\cro{\barmun_{t},\chi}-\cro{\barmu_{t},\chi}\right|\le\xi\Bigl\},
\end{align*}
where $\bar M^{n,\alpha}$ is the martingale defined by (\ref{eq:martX}). 
We have for all $n$, 
\begin{multline}
\label{eq:cvfinal1}
\pr{\sup_{t\in\left(\t2a,T\right]}\left|\cro{\barmun_t,\phi}-\cro{\barmu_t,\phi}\right|>\epsilon}\\
\le \pr{\sup_{t\in\left(\t2a,T\right]}\left|\phi(0)\right|\left|\barmun_t(0)-\barmu_t(0)\right|>{\epsilon \over 2}}
+\pr{\sup_{t\in\left(\t2a,T\right]}\parallel \phi\parallel\left|\cro{\barmun_t,\ind_{\N^*}}-\cro{\barmu_t,\ind_{\N^*}}\right|>{\epsilon \over 2}}\\
\shoveleft{\le \pr{\left\{\sup_{t\in\left(\t2a,T\right]}\left|\barmun_t(0)-\barmu_t(0)\right|>{\epsilon \over 2\mid\phi(0)\mid}\right\}\,
\cap \,\Omega^n_{\alpha,\delta}\,\cap\, \Omega^n_{\alpha,\eta}\, \cap\, \Omega^n_{\alpha,\xi} }}\\
\shoveright{+\pr{\left(\Omega^n_{\alpha,\delta}\right)^c}+\pr{\left(\Omega^n_{\alpha,\eta}\right)^c}+2\pr{\left(\Omega^n_{\alpha,\xi}\right)^c}}\\
+\pr{\left\{\sup_{t\in\left(\t2a,T\right]} \left|\cro{\barmun_t,\ind_{\N^*}}-\cro{\barmu_t,\ind_{\N^*}}\right|>{\epsilon \over 2\parallel\phi\parallel}\right\}\,\cap\,\Omega^n_{\alpha,\xi}}.
\end{multline}
Clearly, from (\ref{eq:afterhorizon2}) and (\ref{eq:condalpha1}), for all $n\in\N^*$, 
\begin{equation}
\label{eq:cvfinal2}
\pr{\left\{\sup_{t\in\left(\t2a,T\right]}\left|\barmun_t(0)-\barmu_t(0)\right|>{\epsilon \over 2\mid\phi(0)\mid}\right\}\,
\cap \,\Omega^n_{\alpha,\delta}\,\cap\, \Omega^n_{\alpha,\eta}\, \cap\, \Omega^n_{\alpha,\xi} }=0.
\end{equation}
On another hand, applying (\ref{eq:convM}) and Proposition \ref{prop:horizonbefore} respectively to $\phi=\ind_0$ and $\phi=\chi$ yields that
\begin{equation}
\label{eq:cvfinal3}
\pr{\left(\Omega^n_{\alpha,\delta}\right)^c}+\pr{\left(\Omega^n_{\alpha,\eta}\right)^c}+2\pr{\left(\Omega^n_{\alpha,\xi}\right)^c}\tend 0.
\end{equation}
Finally, notice that for all $t\ge \t2a$ and for all $n$,
\begin{equation*}
\left|\cro{\barmun_t,\ind_{\N^*}}-\cro{\barmu_t,\ind_{\N^*}}\right|
\le \cro{\barmun_{\t2a},\ind_{\N^*}}+\cro{\barmu_{\t2a},\ind_{\N^*}}
\le \cro{\barmun_{\t2a},\chi}+\cro{\barmu_{\t2a},\chi}\le \cro{\barmun_{\t2a},\chi}+2\alpha,
\end{equation*}
and therefore with (\ref{eq:condalpha2}),
$$
\pr{\left\{\sup_{t\in\left(\t2a,T\right]} \left|\cro{\barmun_t,\ind_{\N^*}}-\cro{\barmu_t,\ind_{\N^*}}\right|>{\epsilon \over 2\parallel\phi\parallel}\right\}\,\cap\,\Omega^n_{\alpha,\xi}}=0.
$$
This together with (\ref{eq:cvfinal2}) and (\ref{eq:cvfinal3}) in (\ref{eq:cvfinal1}), concludes the proof for all $\phi \in \mathcal B_b$ 
such that $\phi(0) \ne 0$.



Only the case where $\phi(0)=0$ remains to be treated. 
For this, we fix $\varepsilon>0$ and let $\alpha>0$ small enough in (\ref{eq:defalpha}), and $\xi>0$ small enough so that
\begin{equation*}
4\alpha+\xi <\varepsilon.
\end{equation*}
Here again, if $T\le \t2a$ then Proposition \ref{prop:horizonbefore} implies that
\begin{equation*}
\pr{\sup_{t\in[0,T]}\left|\cro{\barmun_t,\phi} -\cro{\barmu_t,\phi}\right| >\epsilon}\tend 0.
\end{equation*}
If $T>\t2a$, just write that for all $n$, 
\begin{multline}
\label{eq:cvfinal1phi0}
\pr{\sup_{t\in\left(\t2a,T\right]}\left|\cro{\barmun_t,\phi}-\cro{\barmu_t,\phi}\right|>\epsilon}\\
\le \pr{\left\{\sup_{t\in\left(\t2a,T\right]} \left|\cro{\barmun_t,\phi}-\cro{\barmu_t,\phi}\right|>\epsilon \right\}\,\cap\,\Omega^n_{\alpha,\xi}}
+\pr{\left(\Omega^n_{\alpha,\xi}\right)^c}.
\end{multline}
But then, we have for all $t\ge \t2a$ that
\begin{equation*}
\left|\cro{\barmun_t,\phi}-\cro{\barmu_t,\phi}\right|
\le \parallel \phi \parallel\left(\cro{\barmun_{\t2a},\chi}+2\alpha\right),
\end{equation*}
and hence 
$$\pr{\left\{\sup_{t\in\left(\t2a,T\right]}\left|\cro{\barmun_t,\phi}-\cro{\barmu_t,\phi}\right|>\epsilon \right\}\,\cap\,\Omega^n_{\alpha,\xi}}=0$$
which, together with (\ref{eq:cvfinal3}) in (\ref{eq:cvfinal1phi0}), concludes the proof.
\ep

\subsection{Proof of Corollary \ref{prop:jamming}}
\label{sec:proofjamming}
We conclude with the proof of Corollary \ref{prop:jamming}. Recall the definitions (\ref{eq:defJn}) and (\ref{eq:jammingconstant}).
Let for all $t\ge 0$,
$$c_{\zeta}^t= \lambda\int_0^t \cro{\bar\mu_s,\mathbf 1} ds.$$
Let $\epsilon>0$.
Using simple manipulations of the limiting differential system, we have that for all $t\ge 0$, $\cro{\bar \mu_t,\mathbf 1} \le \exp(-\lambda t)$.
Similarly, applying (\ref{eq:mart}) to $\phi\equiv \mathbf 1$ and taking expectations yields that for any $t\ge 0$ and $n \in \mathbb N^*$,
$$ {d \over dt} \esp{\cro{\mu^n_t,\mathbf 1}} \le - \lambda \esp{\cro{\mu^n_t,\mathbf 1}}.$$
Consequently, there exists $S>0$ such that for all $n\in\N^*$, 
\be \label{ineq:base}
\max\left\{\int_{S}^\infty \cro{\bar \mu_t,\mathbf 1} dt; \int_{S}^\infty \esp{ \cro{\bar\mu^n_t,\mathbf 1}} dt\right\}\le {\epsilon \over 4\lambda}.
\ee
Observe now that for all $n$, $\left(\mu^n,\,J^n\right)$ is a Markov jump process on $\mathcal M_F(\N)\times \N$, whose infinitesimal generator can be
readily deduced from (\ref{eq:gener1}). Applying Dynkin's lemma to the test function
\[F:\left\{\begin{array}{ll}
\mathcal M_F(\N)\times \R& \to \R\\
(\mu,x)&\mapsto x
\end{array}\right.\]
clearly entails that for all $n\in\N^*$ and $t\ge 0$,
$${J^n_t \over n}=\lambda\int_0^t \cro{\bar \mu^n_s,\mathbf 1} ds + {N^n_t \over n},$$
where $N^n$ is a $\mathcal F^n_t$-martingale such that, uniformly over compact time sets
$$ {N^n \over n} \tend \mathbf 0 \text{ in }L^2,$$
as can be proven using Doob's inequality.
We obtain that for $n$ large enough,
\begin{eqnarray*}
\esp{\left|\bar J^n  - c_{\zeta}\right|}&\le& \esp{\left|\frac{J^n_S}{n}  - c^S_{\zeta}\right|}+
\lambda\int_{S}^\infty \cro{\bar\mu_t,\mathbf 1} dt + \lambda\esp{\int_{S}^\infty \cro{\bar\mu^n_t, \mathbf 1} \,dt}\\
&\le& \esp{\left|\lambda\int_0^S\cro{\bar \mu^n_t,\mathbf 1}\,dt  - c^S_{\zeta}\right|}+ \left(\esp{\left({N^n_S \over n}\right)^2}\right)^{1/2} +  \epsilon/2\\
&\le& \epsilon.
\end{eqnarray*}
In the latter, the second inequality follows from (\ref{ineq:base}), the third one from the martingale convergence and
from Theorem \ref{theo.main} applied to $\phi \equiv \mathbf 1$, and by observing that the family of
r.v.'s $\left\{\displaystyle\int_0^S\cro{\bar \mu^n_t,\mathbf 1}\,dt;\,n\in\N^*\right\}$ is bounded by $S$, and hence uniformly integrable.
\ep

\bibliographystyle{plain}
\bibliography{probab2}

\newpage

\appendix

\section{Combinatorial results}
\label{sec:combi}
Let us fix the probability space $\left(\Omega, \F,\mathbf P\right)$. We first introduce the definition, and several basic properties which are used in section \ref{subsec:approx}, for the so called
{\em hypergeometrical distribution}:
\begin{definition}
Let $n, N$ and $p$ be three integers such that $p\ge 1$ and $N\ge 1$.
Let $P:=\left(P(1),...,P(p)\right) \in \N^{p}$ such that $\sum_{i=1}^{p} P(i)=N$. We say that the
measure-valued random variable $Y\in \M_F\left(\N\right)$ follows a {\em multivariate hypergeometrical} distribution of parameters
$\left(n,N,p,P\right)$ if for all $y \in \M_F\left(\N\right)$ of support in $\llbracket 1,p \rrbracket$ such that $y(i)\le P(i)$ for all $i$ 
and $\sum\limits_{i=1}^{p} y(i)=n$, 
$$\pr{Y=y}=\frac{\prod_{i=1}^{p}{P(i) \choose y(i)}}{{N \choose n}}.$$
\end{definition}
The following main characteristics are well-known and easily calculated:
\begin{align}
\esp{Y(i)}&={{n} P(i) \over N},\quad i\in \N\cap [1,n-1]\label{eq:hypergeo1}\\
\mbox{{\bf Cov}}\left(Y(i),Y(j)\right)&={{n}P(i)P(j)\over N^2}{N-{n} \over N-1}, \quad i,j\in \N\cap [1,n-1],\nonumber\\
\esp{Y(i)^3}&={n(n-1)(n-2)}\frac{P(i)^3}{N^3}+3{n(n-1)}\frac{P(i)^2}{N^2}+\frac{{n}P(i)}{N},\quad i\in \N\cap [1,n-1]\label{eq:hypergeo3}.
\end{align}

In particular, we readily deduce from the latter that
\begin{align}
\esp{\cro{Y,\chi}^2}&=\sum_i\sum_j ij\esp{y(i)y(j)}\nonumber\\
&=\sum_i\sum_j ij\Biggl\{\mbox{{\bf Cov}}\left(y(i),y(j)\right)-\esp{y(i)}\esp{y(j)}\Biggl\}\nonumber\\
&=\sum_i\sum_j ij\left\{{nP(i)P(j) \over N^2}{N-{n} \over N-1}-{n^2P(i)P(j) \over N^2}\right\}\nonumber\\
&=\sum_i\sum_j ij{n(n-1)P(i)P(j) \over N(N-1)}.\label{eq:hypergeo2}
\end{align}

\noindent On another hand, as there are at most $p$ integers $i$ such that $y(i)>0$, a simple computation gives
\begin{equation*}
\esp{\cro{Y,\chi}^3} =\sum_{(i,j,\ell)\in [1,n-1]\cap \N}ij\ell \esp{y(i)y(j)y(\ell)} \leq 3{n^2} \sum_{i}i^3 \esp{Y(i)^3}.
\end{equation*}
So with (\ref{eq:hypergeo3}), we obtain
\begin{equation}
\esp{\cro{Y,\chi}^3 } \leq \frac{3{n}^5}{N^3}\sum_{i\in[1,n-1]\cap \N}i^3P(i)^3 +  \frac{3{n}^4}{N^2}\sum_{i\in[1,n-1]}i^3P(i)^2 + \frac{3{n}^3}{N}\sum_{i\in[1,n-1]}i^3P(i).
\label{eq:hypergeo3bis}
\end{equation}

\end{document}

%% file: generator_approximations3.tex
\subsection{Generator approximations}
\label{subsec:approx}
Recall the definition of the map $\Psi:\C\left(\R+,\M_F(\N)\right)\to \C\left(\R+,\M_F(\N)\right)$ in (\ref{eq:limflu}).
We show that the finite variation part ${1\over n}  \mathscr Q^n\Pi_\phi\left(\mu^n\right)$ of $\{\barmun\}$ can be approximated
by $\Psi\left(\barmun\right)$ as long as $\barmun$ takes values in $\M_{\alpha,M}$. To show this, we first write for all $n$, $t$ and $\phi$,
\begin{equation}
\label{eq:decompPsi}
\Psi\left(\barmun\right)_t=:-\lambda\Biggl\{\cro{\barmun_t,\phi}+\cro{\bar A^n\left(\mu^n_t\right),\phi}+\cro{\bar B^n\left(\mu^n_t\right),\phi}\Biggl\},
\end{equation}
where for all $\mu$ such that $\cro{\mu,\chi}>0$,
\begin{align}
\cro{\bar A^n\left(\mu\right),\phi}&={1 \over n}\cro{\mu,\chi\phi};\label{eq:defA*}\\
\cro{\bar B^n\left(\mu\right),\phi}&={1 \over n}\cro{\mu,\chi\Delta\phi}\left(\frac{\cro{\mu,\chi^2}}{\cro{\mu,\chi}}-1\right).
\label{eq:defB*}
\end{align}
Then, we decompose the rescaled generator into two terms $\bar \A^n$ and $\bar \B^n$ that will be shown to get close to $\bar A^n$ and $\bar B^n$ as
$n$ goes large, when applied to a measure $\mu\in n\M_{\alpha,M}$. For doing so, the combinatorial approximation arguments for large $n$ will turn out to be valid only if the first moment of the measure
remains `of order $n$' after the transformation corresponding to the deletion of the new $\sA$-vertex and the new $\sB$-vertices. 
We thus have to tease apart the cases where the degrees of the latter vertices are too big, which correspond to a third term $\bar{\mathscr C}^n$ that will be proven to vanish for large $n$.  
More precisely, we write for all $n$, $\mu \in \M_F(\N)$ and all $\phi \in\mathcal B_b$,
\begin{equation}
\label{eq:decompgenerator}
{1\over n}\mathscr Q^n\Pi_\phi(\mu)=:-\lambda\Biggl\{{1 \over n}\cro{\mu,\phi}+\cro{\bar {\A}^n(\mu),\phi}+\cro{\bar {\B}^n(\mu),\phi}+\cro{\bar {\maC}^n(\mu),\phi}\Biggl\},
\end{equation}
where
\begin{multline}
\cro{\bar \A^n\left(\mu\right),\phi}={1 \over n}\sum_{\substack{k\in\N^*;\\k\le n\alpha/2}}\mu(k)\sum_{\tilde k\in \llbracket 0,k\rrbracket}\P_{\tilde K} (\tilde k| \mu,k)\sum_{\substack{y \in \M_F(\N);\\\cro{y,\mathbf 1}=k}}
\P_Y(y \mid \mu,k,\tilde k)\\
\times\sum_{\tilde y \in \M_F(\N)}
\P_{\tilde Y}(\tilde y \mid \mu,k,\tilde k, y)\cro{\tilde y,\phi};\label{eq:defAn}
\end{multline}
\begin{multline}
\cro{\bar \B^n\left(\mu\right),\phi}={1 \over n}\sum_{\substack{k\in\N^*;\\k\le \cro{\mu,\chi}/4}}\mu(k)\sum_{\tilde k\in \llbracket 0,k\rrbracket}\P_{\tilde K} (\tilde k| \mu,k)\sum_{\substack{y \in \M_F(\N);\\\cro{y,\mathbf 1}=k;\\\cro{y,\chi}\le n\alpha/2}}\P_Y(y \mid \mu,k,\tilde k)\sum_{\tilde y \in \M_F(\N)}
\P_{\tilde Y}(\tilde y \mid \mu,k,\tilde k,y)\\
\times\sum_{i=1}^{n-1}\sum_{x \in \N^*}\P_X(x \mid \mu,k,\tilde k,y,\tilde y)\sum_{w\in\N^{\llbracket 1,n-1 \rrbracket^2}}\P_W(w \mid \mu, k,\tilde k,y,\tilde y, x)
\sum_{\ell=1}^{n-1}\left(\phi(i)-\phi(i-w(i,\ell))\right);\label{eq:defBn}
\end{multline}


\begin{multline}
\cro{\bar {\maC}^n\left(\mu\right),\phi}
={1 \over n}\left(\sum_{\substack{k\in\N^*;\\k> n\alpha/2}}\mu(k)\sum_{\tilde k\in \llbracket 0,k\rrbracket} \P_{\tilde K} (\tilde k| \mu,k)\sum_{\substack{y \in \M_F(\N);\\\cro{y,\mathbf 1}=\tilde k}}\P_Y(y \mid \mu,k,\tilde k)\right.\\
\shoveright{\left.\times\sum_{\tilde y \in \M_F(\N)}\P_{\tilde Y}(\tilde y \mid \mu,k,\tilde k, y)\cro{\tilde y,\phi}\right)}\\
\shoveleft{+{1 \over n}\left(\sum_{\substack{k\in\N^*;\\k\le \cro{\mu,\chi}/4}}\mu(k)\sum_{\tilde k\in \llbracket 0,k\rrbracket} \P_{\tilde K} (\tilde k| \mu,k)\sum_{\substack{y \in \M_F(\N);\\\cro{y,\mathbf 1}=\tilde k;\\\cro{y,\chi}> n\alpha/2}}\P_Y(y \mid \mu,k,\tilde k)\sum_{\tilde y \in \M_F(\N)}
\P_{\tilde Y}(\tilde y \mid \mu,k,k, y)\right.}\\
\shoveright{\left.\times\sum_{i=1}^{n-1}\sum_{x \in \N^*}\P_X(x \mid \mu,k,\tilde k,y,\tilde y)\sum_{w\in\N^{\llbracket 1,n-1 \rrbracket^2}}\P_W(w \mid \mu, k,\tilde k,  y,\tilde y, x)\sum_{\ell=1}^{n-1}\left(\phi(i)-\phi(i-w(i,\ell))\right)\right)}\\
\shoveleft{+{1 \over n}\left(\sum_{\substack{k\in\N^*;\\k> \cro{\mu,\chi}/4}}\mu(k)\sum_{\tilde k\in \llbracket 0,k\rrbracket} \P_{\tilde K} (\tilde k| \mu,k)\sum_{\substack{y \in \M_F(\N);\\\cro{y,\mathbf 1}=\tilde k}}
\P_Y(y \mid \mu,k,\tilde k)\sum_{\tilde y \in \M_F(\N)}
\P_{\tilde Y}(\tilde y \mid \mu,k,k,y)\right.}\\
\shoveright{\left.\times\sum_{i=1}^{n-1}\sum_{x \in \N^*}\P_X(x \mid \mu,k,\tilde k, y, \tilde y)\sum_{w\in\N^{\llbracket 1,n-1 \rrbracket^2}}\P_W(w \mid \mu, k, \tilde k, y,\tilde y, x)\sum_{\ell=1}^{n-1}\left(\phi(i)-\phi(i-w(i,\ell))\right)\right)}\\
=:\cro{\bar \maC^n_1\left(\mu\right),\phi}+\cro{\bar \maC^n_2\left(\mu\right),\phi}+\cro{\bar \maC^n_3\left(\mu\right),\phi}.\label{eq:defEn}
\end{multline}


\noindent We first have the following.
\begin{lemma}
\label{lemma:convCn}
For all sufficiently large $n$, all $\mu\in n\M_{\alpha,M}$ and all bounded $\phi$, 
$$\cro{\bar \maC^n\left(\mu\right),\phi}=o(n).$$
\end{lemma}
\bp
Fix $\phi$, $n>{4 \over 3\alpha}$ and $\mu\in n\M_{\alpha,M}$ throughout the proof. 
Firs, observe that
\begin{equation}
\cro{\bar \maC^n_1(\mu),\phi}
\le {\parallel \phi \parallel \over n}\sum_{\substack{k\in\N^*;\\k> n\alpha/2}}\mu(k)k
\le {\parallel \phi \parallel \over n}{2 \over n\alpha}\sum_{\substack{k\in\N^*;\\k> n\alpha/2}}\mu(k)k^2
\le {2\parallel \phi \parallel \over n^2\alpha}\cro{\mu,\chi^2}
\le {2\parallel \phi \parallel M \over n\alpha}\label{eq:convCn1}.
\end{equation}
\noindent
Second, as $\cro{\mu,\chi}>n\alpha$ we have that
\begin{equation}
\label{eq:oggy1}
k \le {\cro{\mu,\chi} \over 4} \Rightarrow \cro{\mu,\chi}-k \ge {3n\alpha \over 4}.
\end{equation}
Hence, in view of the distribution of $Y(.)$ in (\ref{eq:distriby}) and from (\ref{eq:hypergeo2}), we obtain that for all $k$ and $\tilde k$, 
\begin{equation}
\mathbb E_Y\left[\cro{Y(\mu),\chi}^2\mid \mu,k,\tilde k\right]
\le {k(k-1)\cro{\mu,\chi}^2 \over \left(\cro{\mu,\chi}-k\right)\left(\cro{\mu,\chi} -(k+1)\right)}
\le k(k-1){16M^2 \over 3\alpha\left(3\alpha -4/n\right)}.\label{eq:esperanceL2}
\end{equation}
Therefore, 
\begin{align}
\cro{\bar \maC^n_2(\mu),\phi}
&\le {2\over n}\parallel \phi \parallel\sum_{\substack{k\in\N^*;\\k\le \cro{\mu,\chi}/4}}\mu(k)\sum_{\tilde k\in \llbracket 0,k\rrbracket}\P_{\tilde K} (\tilde k| \mu,k)\sum_{\substack{y \in \M_F(\N);\\\cro{y,\mathbf 1}=\tilde k;\\\cro{y,\chi}> n\alpha/2}}\P_Y(y \mid \mu,k)<y,\chi>\nonumber\\
&\le {2\over n}\parallel \phi \parallel   \sum_{\substack{k\in\N^*;\\k\le \cro{\mu,\chi}/4}}
\mu(k){2 \over n\alpha}\mathbb E_Y\left[\cro{Y(\mu),\chi}^2 \mid \mu,k,\tilde k\right]\nonumber\\
&\le {64\parallel \phi \parallel M^2 \over 3n\alpha^2(3\alpha-4/n)}{1 \over n}\sum_{\substack{k\in\N^*;\\k\le \cro{\mu,\chi}/4}} \mu(k) k(k-1)
\le {64\parallel \phi \parallel M^3 \over 3n\alpha^2(3\alpha-4/n)},\label{eq:convCn2}
\end{align}
where we reason as for (\ref{eq:convCn1}) in the second upper-bound.

\noindent
Last, as for all $k$ and $y$, $\cro{\mu,\chi}-k-\cro{y,\chi} \ge 0$, we have that
\begin{equation}
\label{eq:oggy2}
k > \cro{\mu,\chi}/4 \Rightarrow \cro{y,\chi}\le {3\cro{\mu,\chi} \over 4} <3k.
\end{equation}
Therefore, 
\begin{align}
\cro{\bar \maC^n_3(\mu),\phi}
&\le 2\parallel \phi \parallel {1 \over n}\sum_{\substack{k\in\N^*;\\k> \cro{\mu,\chi}/4}}\mu(k)\sum_{\tilde k\in \llbracket 0,k\rrbracket}\P_{\tilde K} (\tilde k| \mu,k)\sum_{\substack{y \in \M_F(\N);\\\cro{y,\mathbf 1}=\tilde k}}\P_Y(y \mid \mu,k,\tilde k)\cro{y,\chi}\nonumber\\
&\le 2\parallel \phi \parallel {1 \over n}\sum_{\substack{k\in\N^*;\\k> \cro{\mu,\chi}/4}}
\mu(k)3k \le 2\parallel \phi \parallel {1 \over n}\sum_{\substack{k\in\N^*;\\k> n\alpha/4}}
\mu(k)3k\le {24\parallel \phi \parallel M \over n\alpha},\label{eq:convCn3}
\end{align}
applying again the same argument as for (\ref{eq:convCn1}). The proof is complete by gathering (\ref{eq:convCn1}), (\ref{eq:convCn2}) and (\ref{eq:convCn3}).
\ep

\noindent We now turn to the terms $\bar \A^n$ and  $\bar \B^n$.
\begin{lemma}
\label{lemma:convAn}
For all sufficient large $n$, all $\mu \in n\M_{\alpha,M}$ and all bounded $\phi$,
$$\left|\cro{\bar \A^n(\mu),\phi}-\cro{\bar A^n(\mu),\phi}\right|=o(n).$$
\end{lemma}
\bp
Fix $\phi$, $n>{2\over\alpha}$ and $\mu \in n\M_{\alpha,M}$ throughout the proof. 
First, for all $k$ such that $k \le n\alpha/2$ we have that
$$\cro{\mu,\chi} -k \ge {n\alpha \over 2},$$
which will be used in several steps of this proof.
In particular, from (\ref{eq:majorePytilde}) we clearly have that
\begin{equation}
\mathbb P_{\tilde Y}(y |\mu,k,k,y)\ge 1-{k \choose 2}4\frac{M}{n\alpha\left(\alpha-{2\over n}\right)}\cdot\label{eq:convAn0}
\end{equation}
\noindent
 Recall (\ref{eq:Pkktilde}), and let us rewrite 
\begin{multline}
\label{eq:convA0}
\cro{\bar \A^n\left(\mu\right),\phi}={1 \over n}\sum_{\substack{k\in\N^*;\\k\le n\alpha/2}}\mu(k)\P_{\tilde K}(k \mid \mu,k)
\sum_{\substack{y \in \M_F(\N);\\\cro{y,\mathbf 1}=k}}\P_Y(y \mid \mu,k,k)\sum_{\tilde y \in \M_F(\N)}\P_{\tilde Y}(\tilde y \mid \mu,k,k, y)\cro{\tilde y,\phi}\\
+{1 \over n}\sum_{\substack{k\in\N^*;\\k\le n\alpha/2}}\mu(k)\sum_{\tilde k\in \llbracket 0,k-1\rrbracket} \P_{\tilde K} (\tilde k| \mu,k)\sum_{\substack{y \in \M_F(\N);\\\cro{y,\mathbf 1}=\tilde k}}\P_Y(y \mid \mu,k,\tilde k)
\sum_{\tilde y \in \M_F(\N)}\P_{\tilde Y}(\tilde y \mid \mu,k,\tilde k, y)\cro{\tilde y,\phi}\\
\shoveleft{={1 \over n}\sum_{\substack{k\in\N^*;\\k\le n\alpha/2}}\mu(k)
\sum_{\substack{y \in \M_F(\N);\\\cro{y,\mathbf 1}=k}}\P_Y(y \mid \mu,k,k)\sum_{\tilde y \in \M_F(\N)}\P_{\tilde Y}(\tilde y \mid \mu,k,k, y)\cro{\tilde y,\phi}}\\
+{1 \over n}\sum_{\substack{k\in\N^*;\\k\le n\alpha/2}}\mu(k)\left(\P_{\tilde K}(k \mid \mu,k)-1\right)
\sum_{\substack{y \in \M_F(\N);\\\cro{y,\mathbf 1}=k}}\P_Y(y \mid \mu,k,k)\sum_{\tilde y \in \M_F(\N)}\P_{\tilde Y}(\tilde y \mid \mu,k,k, y)\cro{\tilde y,\phi}\\
+{1 \over n}\sum_{\substack{k\in\N^*;\\k\le n\alpha/2}}\mu(k)\sum_{\tilde k\in \llbracket 0,k-1\rrbracket} \P_{\tilde K} (\tilde k| \mu,k)\sum_{\substack{y \in \M_F(\N);\\\cro{y,\mathbf 1}=\tilde k}}\P_Y(y \mid \mu,k,\tilde k)
\sum_{\tilde y \in \M_F(\N)}\P_{\tilde Y}(\tilde y \mid \mu,k,\tilde k, y)\cro{\tilde y,\phi}\\
=:\cro{\bar \A^n_1(\mu),\phi}+\cro{\bar \A^n_2(\mu),\phi}+\cro{\bar \A^n_3(\mu),\phi}.
\end{multline}
The last two terms vanish as a consequence of the fact that $\tilde K(\mu)$ and $K(\mu)$ tend to coincide as the size of the graph goes large. To see this, 
from (\ref{eq:Pkktilde}) we have that for all $k \le n\alpha/2$, 
\begin{align}\P_{\tilde K}(k \mid \mu,k)&
\ge {{\cro{\mu,\chi}-k  \choose k} \over {\cro{\mu,\chi} \choose k}}\ind_{\{k \le \cro{\mu,\chi}-k\}}\nonumber\\
&={{\cro{\mu,\chi}-k  \choose k} \over {\cro{\mu,\chi} \choose k}}\nonumber\\
&={\left(\cro{\mu,\chi}-k\right)....\left(\cro{\mu,\chi}-2k+1\right) \over \left(\cro{\mu,\chi}\right)....\left(\cro{\mu,\chi}-k+1\right)}\nonumber\\
&\ge 
\left(1 - {k \over \cro{\mu,\chi}-k+1}\right)^k\nonumber\\
&\ge 1-{k^2 \over \cro{\mu,\chi}-k+1}.
\label{eq:convKtildeK}
\end{align}
\noindent Therefore, 
as $\cro{\tilde y,\mathbf 1}\le k$ for any $k,\tilde k$ and $\tilde y$ we have that
\begin{equation}
\left|\cro{\bar \A^n_2(\mu),\phi}+\cro{\bar \A^n_3(\mu),\phi}\right|\le 
{2 \over n}\sum_{\substack{k\in\N^*;\\k\le n\alpha/2}}\mu(k)k\left(1-\P_{\tilde K}(k \mid \mu,k)\right)
\le {2 \over n}{\cro{\mu,\chi^3} \over \cro{\mu,\chi}-k+1}\le {4M \over n\alpha}.\label{eq:convA23}
\end{equation}

\bigskip

\noindent
Let us now address the term $\bar \A^n_1$. We have
\begin{multline}
\label{eq:convAn1}
\cro{\bar \A^n_1(\mu),\phi}
={1\over n}\sum_{\substack{k\in\N^*;\\k\le n\alpha/2}}\mu(k)\sum_{\substack{y \in \M_F(\N);\\\cro{y,\mathbf 1}=k}}\P_Y(y \mid \mu,k)
\mathbb E_{\tilde Y}\left[\cro{\tilde Y(\mu),\phi}\ind_{\{\tilde Y(\mu) \ne y\}} \mid \mu,k,y\right]\\
\shoveright{+{1\over n}\sum_{\substack{k\in\N^*;\\k\le n\alpha/2}}\mu(k)\sum_{y \in \M_F(\N)}\P_Y(y \mid \mu,k)\P_{\tilde Y}\left(y \mid \mu,k,k,y\right)
\cro{y,\phi}}\\
\shoveleft{={1\over n}\left\{\sum_{\substack{k\in\N^*;\\k\le n\alpha/2}}\mu(k)\sum_{\substack{y \in \M_F(\N);\\\cro{y,\mathbf 1}=k}}\P_Y(y \mid \mu,k)\right.}\\
\shoveright{\times\Biggl(\mathbb E_{\tilde Y}\left[\cro{\tilde Y(\mu),\phi}\ind_{\{\tilde Y(\mu) \ne y\}} \mid \mu,k,y\right]+\cro{y,\phi}\Bigl(\P_{\tilde Y}\left(y \mid \mu,k,k,y\right)-1\Bigl)\Biggl)\Biggl\}}\\
\shoveright{+{1\over n}\sum_{\substack{k\in\N^*;\\k\le n\alpha/2}}\mu(k)\mathbb E_Y\left[\cro{Y(\mu),\phi} \mid \mu,k\right]}\\
=:\cro{\bar \A^n_{11}(\mu),\phi}+\cro{\bar \A^n_{12}(\mu),\phi}. 
\end{multline}
\noindent
First, as $\cro{y,\mathbf 1}\le k$, from (\ref{eq:convAn0}) we have that 
\begin{align}
\left|\cro{\bar \A^n_{11}(\mu),\phi}\right|
&\le  {2\parallel \phi \parallel\over n}\sum_{\substack{k\in\N^*;\\k\le n\alpha/2}}\mu(k)\sum_{y \in \M_F(\N)}\P_Y(y \mid \mu,k)
\Bigl(1-\P_{\tilde Y}\left(y \mid \mu,k,k,y\right)\Bigl)k\nonumber\\
&\le  {2\parallel \phi \parallel\over n}\sum_{\substack{k\in\N^*;\\k\le n\alpha/2}}\mu(k){k \choose 2}4\frac{M}{n\alpha\left(\alpha-{2 \over n}\right)}k\le \frac{4\parallel \phi \parallel M^2}{n\alpha\left(\alpha-{2 \over n}\right)}.\label{eq:convAn2}
\end{align}

\noindent
Second, using again the hypergeometric distribution of $Y(.)$ and (\ref{eq:hypergeo1}) we have that
\begin{multline*}
\begin{aligned}
\left|\cro{\bar \A^n_{12}(\mu),\phi}-\cro{\bar A^n(\mu),\phi}\right|
&={1\over n}\left|\sum_{\substack{k\in\N^*;\\k\le n\alpha/2}}\mu(k)\sum_{i=1}^{n-1}\phi(i)\mathbb E_Y\left[Y(\mu)(i)\mid \mu,k\right]-\cro{\mu,\phi\chi}\right|\\
&={1\over n}\left|\sum_{\substack{k\in\N^*;\\k\le n\alpha/2}}\mu(k)\sum_{i=1}^{n-1}\phi(i)k{i\mu(i)-i\delta_k(i) \over \cro{\mu,\chi}-k}-\cro{\mu,\phi\chi}\right|
\\
&\le {1\over n}\sum_{\substack{k\in\N^*;\\k\le n\alpha/2}}\mu(k){\parallel \phi \parallel k \over \cro{\mu,\chi} -k}
+{\cro{\mu,\phi\chi}\over n}\left|\sum_{\substack{k\in\N^*;\\k\le n\alpha/2}}{\mu(k)k\over \cro{\mu,\chi}-k}-1\right|
\end{aligned}\\
\shoveright{\le {2M\parallel \phi \parallel \over n\alpha}
+ {\cro{\mu,\phi\chi}\over n}\sum_{\substack{k\in\N^*;\\k\le n\alpha/2}}{\mu(k)k^2 \over \cro{\mu,\chi}\left(\cro{\mu,\chi} -k\right)}
+{\cro{\mu,\phi\chi}\over n}\sum_{\substack{k\in\N^*;\\ k> n\alpha/2}}{\mu(k)k \over \cro{\mu,\chi}}}\\
\le {2M\parallel \phi \parallel \over n\alpha}+{2M^2\parallel \phi \parallel \over n\alpha^2}
+{M\parallel \phi \parallel \over \alpha}{2M \over n\alpha},
\end{multline*}
where the upper-bound of the third term of the last inequality follows from the same argument as for (\ref{eq:convCn1}).
This, together with (\ref{eq:convAn2}) in (\ref{eq:convAn1}), and then with (\ref{eq:convA23}) in (\ref{eq:convA0}), concludes the proof.
\ep

\begin{lemma}
\label{lemma:convBn}
For all $n$, all $\mu\in n\M_{\alpha,M}$ and all bounded $\phi$,
$$\left|\cro{\bar \B^n(\mu),\phi}-\cro{\bar B^n\left(\mu\right),\phi}\right|=o(n).$$
\end{lemma}

\bp
Fix $\phi\in \mathcal B_b$, $n>{4\over \alpha}$ and $\mu\in n\M_{\alpha,M}$. 
We split $\bar \B^n$ as follows:  
\begin{equation}
\label{eq:defBn0}
\cro{\bar \B^n\left(\mu\right),\phi}=\sum_{i=1}^4 \cro{\bar \B^n_i\left(\mu\right),\phi},
\end{equation}
where 
\begin{multline}
\label{eq:defBn1}
\cro{\bar \B^n_1\left(\mu\right),\phi}={1 \over n}\left(\sum_{\substack{k\in\N^*;\\k\le \cro{\mu,\chi}/4}}\mu(k)\sum_{\substack{y \in \M_F(\N);\\\cro{y,\mathbf 1}=k;\\\cro{y,\chi}\le n\alpha/2}}\P_Y(y \mid \mu,k,k)\sum_{x \in \N^*}\P_X(x \mid \mu,k,k,y,y)\right.\\
\left.\times\sum_{w\in\N^{\llbracket 1,n-1 \rrbracket^2}}\P_W(w \mid \mu, k,k,y,y, x)\sum_{i=1}^{n-1}Z(\mu)(i)\Delta\phi\left(i\right)\right);
\end{multline}

\begin{multline}
\label{eq:defBn2}
\cro{\bar \B^n_2\left(\mu\right),\phi}={1 \over n}\left(\sum_{\substack{k\in\N^*;\\k\le \cro{\mu,\chi}/4}}\mu(k)\sum_{\substack{y \in \M_F(\N);\\\cro{y,\mathbf 1}=k;\\\cro{y,\chi}\le n\alpha/2}}\P_Y(y \mid \mu,k,k)\sum_{x \in \N^*}\P_X(x \mid \mu,k,k,y,y)\right.\\
\left.\times\sum_{\substack{ w\in\N^{\llbracket 1,n-1 \rrbracket^2};\\\exists (i_0,\ell_0) \in \llbracket 1,n-1 \rrbracket^2; w(i,\ell_0)\ge 2}}\P_W(w \mid \mu, k,k,y,y, x)\sum_{i=1}^{n-1}\left(\sum_{\ell=1}^{n-1}\left(\phi(i)-\phi(i-w(i,\ell))\right)-Z(\mu)(i)\Delta\phi\left(i\right)\right)\right);
\end{multline}

\begin{multline}
\label{eq:defBn3}
\cro{\bar \B^n_3\left(\mu\right),\phi}={1 \over n}\left(\sum_{\substack{k\in\N^*;\\k\le \cro{\mu,\chi}/4}}\mu(k)\sum_{\substack{y \in \M_F(\N);\\\cro{y,\mathbf 1}=k;\\\cro{y,\chi}\le n\alpha/2}}\P_Y(y \mid \mu,k,k)
\left(\P_{\tilde Y}(y \mid \mu,k,k,y)-1\right)\right.\\
\left.\times\sum_{x \in \N^*}\P_X(x \mid \mu,k,k,y,y)\sum_{w\in\N^{\llbracket 1,n-1 \rrbracket^2}}\P_W(w \mid \mu, k,k,y,y, x)
\sum_{i=1}^{n-1}\sum_{\ell=1}^{n-1}\left(\phi(i)-\phi(i-w(i,\ell))\right)\right)\\
\shoveleft{+{1 \over n}\left(\sum_{\substack{k\in\N^*;\\k\le \cro{\mu,\chi}/4}}\mu(k)\sum_{\substack{y \in \M_F(\N);\\\cro{y,\mathbf 1}=k;\\\cro{y,\chi}\le n\alpha/2}}\P_Y(y \mid \mu,k,k)\sum_{\substack{\tilde y \in \M_F(\N);\\\tilde y \ne y}}
\P_{\tilde Y}(\tilde y \mid \mu,k,k,y)\right.}\\
\left.\times\sum_{x \in \N^*}\P_X(x \mid \mu,k,k,y,\tilde y)\sum_{w\in\N^{\llbracket 1,n-1 \rrbracket^2}}\P_W(w \mid \mu, k,k,y,\tilde y, x)
\sum_{i=1}^{n-1}\sum_{\ell=1}^{n-1}\left(\phi(i)-\phi(i-w(i,\ell))\right)\right);
\end{multline}

\begin{multline}
\label{eq:defBn4}
\cro{\bar \B^n_4\left(\mu\right),\phi}={1 \over n}\left(\sum_{\substack{k\in\N^*;\\k\le \cro{\mu,\chi}/4}}\mu(k)\left(\P_{\tilde K} (k| \mu,k)-1\right)\sum_{\substack{y \in \M_F(\N);\\\cro{y,\mathbf 1}=k;\\\cro{y,\chi}\le n\alpha/2}}\P_Y(y \mid \mu,k,k)\sum_{\tilde y \in \M_F(\N)}
\P_{\tilde Y}(\tilde y \mid \mu,k,\tilde k,y)\right.\\
\left.\times\sum_{x \in \N^*}\P_X(x \mid \mu,k,k,y,\tilde y)\sum_{w\in\N^{\llbracket 1,n-1 \rrbracket^2}}\P_W(w \mid \mu, k,k,y,\tilde y, x)
\sum_{i=1}^{n-1}\sum_{\ell=1}^{n-1}\left(\phi(i)-\phi(i-w(i,\ell))\right)\right)\\
\shoveleft{+{1 \over n}\left(\sum_{\substack{k\in\N^*;\\k\le \cro{\mu,\chi}/4}}\mu(k)\sum_{\tilde k\in \llbracket 0,k-1\rrbracket}\P_{\tilde K} (\tilde k| \mu,k)\sum_{\substack{y \in \M_F(\N);\\\cro{y,\mathbf 1}= \tilde k;\\\cro{y,\chi}\le n\alpha/2}}\P_Y(y \mid \mu,k,\tilde k)\sum_{\tilde y \in \M_F(\N)}
\P_{\tilde Y}(\tilde y \mid \mu,k,\tilde k,y)\right.}\\
\left.\times\sum_{x \in \N^*}\P_X(x \mid \mu,k,\tilde k,y,\tilde y)\sum_{w\in\N^{\llbracket 1,n-1 \rrbracket^2}}\P_W(w \mid \mu, k,\tilde k,y,\tilde y, x)
\sum_{i=1}^{n-1}\sum_{\ell=1}^{n-1}\left(\phi(i)-\phi(i-w(i,\ell))\right)\right).
\end{multline}
In the latter partition, the last three terms $\bar \B^n_4\left(.\right)$, $\bar \B^n_3\left(.\right)$ and $\bar \B^n_2\left(.\right)$, which 
gather respectively the cases of existence of self-loops around the new $\sA$-node, of multiple edges between the new $\sA$-node and the new $\sB$-nodes, 
and of multi-edges between the new $\sB$-nodes and the remaining $\sU$-nodes, will be shown to be small w.r.t. $n$. Only the term $\bar \B^n_1\left(.\right)$ is non degenerate, and as we will prove, is at a distance to $\bar B^n\left(.\right)$ which is also small w.r.t. $n$.  

To prove this, let us 
first recall the definition of $Z(\mu)$ in (\ref{eq:defZ}) and its distribution in (\ref{eq:distribz}): for any $i$, there are at most $Z(\mu)(i)$ 
indexes $\ell \in \llbracket 1,n-1 \rrbracket$ such that $W(\mu)(i,\ell)$ differs from zero. Therefore we have that
\begin{align}
\sum_{i=1}^{n-1}\sum_{\ell=1}^{n-1}\left(\phi(i)-\phi(i-w(i,\ell))\right)
&\le 2\parallel \phi \parallel \sum_{i=1}^{n-1}Z(\mu)(i)\nonumber\\
&\le 2\parallel \phi \parallel X(\mu)
\le 2\parallel \phi \parallel \left(\cro{Y(\mu),\chi}-K(\mu)\right).\label{eq:compareZ}
\end{align} 
On another hand, notice that for all $k$, 
\begin{multline}
\sum_{\substack{y \in \M_F(\N);\\\cro{y,\mathbf 1}=k;\\\cro{y,\chi}\le n\alpha/2}}\P_Y(y \mid \mu,k)\left(\cro{y,\chi}-k\right)\\
\begin{aligned}
&=\mathbb E_Y\left[\cro{Y(\mu),\chi}-k \mid \mu,k\right]-\sum_{\substack{y \in \M_F(\N);\\\cro{y,\mathbf 1}=k;\\\cro{y,\chi}> n\alpha/2}}\P_Y(y \mid \mu,k)\left(\cro{y,\chi}-k\right)\\
&=\sum_{i=1}^{n-1}ik \frac{i\left(\mu(i)-\delta_k(i)\right)}{\cro{\mu,\chi}-k}-k-\sum_{\substack{y \in \M_F(\N);\\\cro{y,\mathbf 1}=k;\\\cro{y,\chi}> n\alpha/2}}\P_Y(y \mid \mu,k)\left(\cro{y,\chi}-k\right)\\
&= k\left({\cro{\mu,\chi^2} \over \cro{\mu,\chi}-k}-1\right)-{k^3 \over \cro{\mu,\chi}-k}-\sum_{\substack{y \in \M_F(\N);\\\cro{y,\mathbf 1}=k;\\\cro{y,\chi}> n\alpha/2}}\P_Y(y \mid \mu,k)\left(\cro{y,\chi}-k\right)
\end{aligned}\\
=k\left({\cro{\mu,\chi^2} \over \cro{\mu,\chi}-k}-1\right)+k^3o_3(n)+k(k-1)o_4(n),
\label{eq:esperanceLbis}
\end{multline}
where the equivalent $o_4(n)$ follows from (\ref{eq:esperanceL2}), just as in (\ref{eq:convCn2}). 
Clearly, (\ref{eq:esperanceLbis}) implies in particular that
\begin{equation}
\sum_{\substack{y \in \M_F(\N);\\\cro{y,\mathbf 1}=k;\\\cro{y,\chi}\le n\alpha/2}}\P_Y(y \mid \mu,k)\left(\cro{y,\chi}-k\right)
\le k{4M \over 3\alpha}.\label{eq:esperanceL}
\end{equation}

We are now in position to address the vanishing terms in the decomposition (\ref{eq:defBn0}), from bottom to top. First, as for (\ref{eq:convA23}), the term $\bar \B^n_{4}(.)$ vanishes since $\tilde K(\mu)$ and $K(\mu)$ tend to coincide. 
More precisely, we have with (\ref{eq:compareZ}) that 
\begin{align}
\left|\cro{\bar \B^n_{4}(\mu),\phi} \right|  
&\leq   {4 \parallel\phi\parallel\over n}\sum_{\substack{k\in\N^*;\\k\le \cro{\mu,\chi}/4}}\mu(k) \left(1-\P_{\tilde K} (k| \mu,k)\right)\sum_{\substack{y \in \M_F(\N);\\\cro{y,\mathbf 1}= k;\\\cro{y,\chi}\le n\alpha/2}}\P_Y(y \mid \mu,k, k) (\cro{y,\chi}-k) \nonumber \nonumber\\
&\leq   {4 \parallel\phi\parallel \over n}\sum_{\substack{k\in\N^*;\\k\le \cro{\mu,\chi}/4}}\mu(k) \frac{k^2}{\cro{\mu,\chi}-k+1}\sum_{\substack{y \in \M_F(\N);\\\cro{y,\mathbf 1}= k;\\\cro{y,\chi}\le n\alpha/2}}\P_Y(y \mid \mu,k, k) (\cro{y,\chi}-k)\nonumber\\
&\leq   {4 \parallel\phi\parallel\over n} \frac{2}{2n\alpha+1}\frac{4M}{3\alpha}\sum_{\substack{k\in\N^*;\\k\le \cro{\mu,\chi}/4}}\mu(k) k^3 \leq 4 \parallel\phi\parallel \frac{2}{2n\alpha+1}\frac{4M^2}{3\alpha},
\label{eq:convBn4}
\end{align}
where we use \eqref{eq:convKtildeK} in the second upper-bound and (\ref{eq:esperanceL}) in the last one.  

The term $\bar \B^n_{3}(.)$ can be treated analogously. Applying again (\ref{eq:compareZ}) and then (\ref{eq:convAn0}), we have that
\begin{align} 
\left|\cro{\bar \B^n_{3}(\mu),\phi} \right|  
&\le   {4 \parallel\phi\parallel\over n}\sum_{\substack{k\in\N^*;\\k\le \cro{\mu,\chi}/4}}\mu(k) \sum_{\substack{y \in \M_F(\N);\\\cro{y,\mathbf 1}= k;\\\cro{y,\chi}\le n\alpha/2}}\P_Y(y \mid \mu,k, k)\left(1- \P_{\tilde Y}(y \mid \mu,k,k,y)\right) (\cro{y,\chi}-k)\nonumber\nonumber\\
& \le   {4 \parallel\phi\parallel\over n}\sum_{\substack{k\in\N^*;\\k\le \cro{\mu,\chi}/4}}\mu(k){k \choose 2}4\frac{M}{n\alpha\left(\alpha-{2\over n}\right)} \sum_{\substack{y \in \M_F(\N);\\\cro{y,\mathbf 1}= k;\\\cro{y,\chi}\le n\alpha/2}}\P_Y(y \mid \mu,k, k)(\cro{y,\chi}-k)\nonumber\\
&\le   {8 \parallel\phi\parallel\over n}\frac{M}{n\alpha\left(\alpha-{2\over n}\right)}{4M \over 3\alpha}
\sum_{\substack{k\in\N^*;\\k\le \cro{\mu,\chi}/4}}\mu(k)k^2
\le {32 \parallel\phi\parallel\over 3n}\frac{M^3}{\alpha^2\left(\alpha-{2\over n}\right)},
\label{eq:convBn3}
\end{align}
where we use again (\ref{eq:esperanceL}) in the third upper-bound.

\noindent 
We now examine the term $\bar\B^n_2(.)$ defined by (\ref{eq:defBn2}). For all
$k\le \cro{\mu,\chi}/4$ and $y$ such that $\cro{y,\chi}\le {n\alpha \over 2}$, we have that
\begin{equation}
\label{eq:TRUC}
\cro{\mu,\chi}-\left(k+\cro{y,\chi}\right)\ge {n\alpha \over 4}.
\end{equation} 
Thus, likewise (\ref{eq:convAn0}), from (\ref{eq:majorecheckQ}) we have that 
\begin{equation*}
\mathbb Q_W\left(\mu,k,k,y,y,x\right)\ge 1-{x \choose 2}{M \over {n\alpha \over 4}\left({\alpha \over 4}-{1\over n}\right)}=: 1-x^2o_1(n).
\label{eq:convBnQ}
\end{equation*}

\noindent Therefore, from the second inequality in (\ref{eq:compareZ}) we have that 
\begin{multline*}
\left|\cro{\bar \B^n_2(\mu),\phi}\right|
\le {2\parallel \phi \parallel \over n}\left(\sum_{\substack{k\in\N^*;\\k\le \cro{\mu,\chi}/4}}\mu(k)
 \sum_{\substack{y \in \M_F(\N);\\\cro{y,\mathbf 1}=k;\\\cro{y,\chi}\le n\alpha/2}}\P_Y(y \mid \mu,k,k)\right.\\
\shoveright{\left.\times\sum_{x \in \N^*}\P_X(x \mid \mu,k,k, y,y)\left(1-\mathbb Q_W\left(\mu,k,k,y,y,x\right)\right) 2x\right)}\\
\leq {4\parallel \phi \parallel \over n} o_1(n) {1\over n}\sum_{\substack{k\in\N^*;\\k\le \cro{\mu,\chi}/4}}\mu(k)\mathbb E_Y\left[\cro{Y(\mu),\chi}^3 \mid \mu,k,k \right].
\end{multline*}

\noindent But \eqref{eq:hypergeo3} implies that
\begin{align*}
\mathbb E_Y\left[\cro{Y(\mu),\chi}^3 \mid \mu,k,k \right]
&\le \frac{3 k^5  \cro{\mu,\mathbf 1}^2 \cro{\mu,\chi^6}}{\left(\cro{\mu,\chi}-k\right)^3}  +  \frac{3 k^4  \cro{\mu,\mathbf 1} \cro{\mu,\chi^5}}
{\left(\cro{\mu,\chi}-k\right)^2} + \frac{3 k^3 \cro{\mu,\chi^4}}{\cro{\mu,\chi}-k} \\
&\le \frac{3k^5M^3 4^3}{\alpha^3} +\frac{3k^4M^24^2}{\alpha^2}+\frac{12k^3M}{\alpha},
\end{align*}
so injecting this in \eqref{eq:convBn2} we obtain that
\begin{equation}
\label{eq:convBn2}
\left|\cro{\bar \B^n_2(\mu),\phi}\right|\le 4\parallel \phi \parallel C o_1(n),
\end{equation}
where the constant $C$ is given by
\begin{gather*}
C=\frac{3	M^4 4^3}{\alpha^3} +\frac{3M^3  4^2}{\alpha^2}+\frac{12M^2}{\alpha}.
\end{gather*}

Let us finally focus on the first term (\ref{eq:defBn1}). Recalling the distribution of $Z(.)$ in (\ref{eq:distribz}) and using 
(\ref{eq:hypergeo1}), the latter can be rewritten as follows, 
\begin{multline}
\label{eq:convBn1}
\cro{\bar \B^n_{1}(\mu),\phi}\\
\shoveleft{=\frac{1}{n}\sum_{\substack{k\in\N^*;\\k\le \cro{\mu,\chi}/4}}\mu(k) \sum_{\substack{y \in \M_F(\N);\\\cro{y,\mathbf 1}=k;\\\cro{y,\chi}\le n\alpha/2}}\P_Y(y \mid \mu,k,k) \mathbb E_X\left[X(\mu) \mid \mu,k,k,y,y \right] }\sum_{i=1}^{n-1}\Delta\phi(i)\frac{ i\left(\mu(i)-\delta_{k}(i)-y(i)\right)}{ \cro{\mu,\chi}-\cro{y,\chi}-k}\\
\shoveleft{=\frac{1}{n}\sum_{\substack{k\in\N^*;\\k\le \cro{\mu,\chi}/4}}\mu(k)
\sum_{\substack{y \in \M_F(\N);\\\cro{y,\mathbf 1}= k;\\\cro{y,\chi}\le n\alpha/2}}\P_Y(y \mid \mu,k, k) (\cro{Y(\mu),\chi}- k)}
\sum_{i=1}^{n-1}\Delta\phi(i)\frac{ i\left(\mu(i)-\delta_{k}(i)-y(i)\right)}{ \cro{\mu,\chi}-\cro{y,\chi}-k}\\
\shoveleft{+\frac{1}{n}\sum_{\substack{k\in\N^*;\\k\le \cro{\mu,\chi}/4}}{\mu(k)} 
\sum_{\substack{y \in \M_F(\N);\\\cro{y,\mathbf 1}= k;\\\cro{y,\chi}\le n\alpha/2}}\P_Y(y \mid \mu,k,k) }\sum_{i=1}^{n-1}\Delta\phi(i)\frac{ i\left(\mu(i)-\delta_{k}(i)-y(i)\right)}{ \cro{\mu,\chi}-\cro{y,\chi}-k}\\
\times\Biggl(\mathbb E_X\left[X(\mu)  \ind_{\{X(\mu)\neq \cro{y,\chi}-k\}}\mid \mu,k, k,y,y \right] 
+ (\cro{y,\chi}- k)\Bigl(\mathbb{P}_{X}\left(\cro{y,\chi}-k |\mu,k,k,y,y\right)-1\Bigl) \Biggl)\\
=:\cro{\bar \B^n_{11}(\mu),\phi}+\cro{\bar \B^n_{12}(\mu),\phi}.
\end{multline}

\noindent Let us prove that $\cro{\bar \B^n_{12}(\mu),\phi}$ also vanishes. Likewise \eqref{eq:convKtildeK}, we have that for all 
$k$ and $y$,  
\begin{gather*}
\mathbb P_{X}\left(\cro{y,\chi}-k |\mu,k,k,y,y\right)\geq 1-\frac{(\cro{y,\chi}-k)^2}{\cro{\mu,\chi}-\cro{y,\chi}-k+1}.
\end{gather*}
This implies using (\ref{eq:TRUC}) that 
\begin{multline}
\cro{\bar \B^n_{12}(\mu),\phi} 
\leq \frac{4\parallel \phi \parallel}{n}\sum_{\substack{k\in\N^*;\\k\le \cro{\mu,\chi}/4}}{\mu(k)} 
\sum_{\substack{y \in \M_F(\N);\\\cro{y,\mathbf 1}= k;\\\cro{y,\chi}\le n\alpha/2}}\P_Y(y \mid \mu,k,k) \frac{(\cro{y,\chi}-k)^2}{\cro{\mu,\chi}-\cro{y,\chi}-k+1} (\cro{y,\chi}-k)\\
\leq \frac{16\parallel \phi \parallel}{n(n\alpha+4)}\sum_{\substack{k\in\N^*;\\k\le \cro{\mu,\chi}/4}}{\mu(k)} 
 \E_Y\left[\cro{Y(\mu),\chi}^3|\mu,k,k\right],
\end{multline}
which is again a $o(n)$, just as \eqref{eq:convBn2}. 
It only remains to prove that $\left|\cro{\bar \B^n_{11}(\mu),\phi}-\cro{\bar B^n(\mu),\phi}\right|$ also is a $o(n)$. 
We have
 \begin{multline}
\label{eq:convBn11}
\left|\cro{\bar \B^n_{11}(\mu),\phi}-\cro{\bar B^n(\mu),\phi}\right|\\ 
\leq \frac{1}{n} \left| \cro{\mu,\chi\Delta \phi} \sum_{\substack{k\in\N^*;\\k\le \cro{\mu,\chi}/4}}{\mu(k)}  \sum_{\substack{y \in \M_F(\N);\\\cro{y,\mathbf 1}= k;\\\cro{y,\chi}\le n\alpha/2}}\P_Y(y \mid \mu,k,k) 
\left({\cro{y,\chi}-k \over\cro{\mu,\chi}-\cro{y,\chi}-k}  \right) -\left({\cro{\mu,\chi^2} \over \cro{\mu,\chi}}-1\right) \right|\\
+ \frac{1}{n} \left| \sum_{\substack{k\in\N^*;\\k\le \cro{\mu,\chi}/4}}{\mu(k)}\sum_{\substack{y \in \M_F(\N);\,\cro{y,\mathbf 1}= k;\\\cro{y,\chi}\le n\alpha/2}}\P_Y(y \mid \mu,k, k) 
\left(\cro{y,\chi}- k \right)\sum_{i=1}^{n-1}\Delta\phi(i)\frac{ i\left(\delta_{k}(i)+y(i)\right)}{ \cro{\mu,\chi}-\cro{y,\chi}-k} 
\right|\\
\shoveleft{\leq  \frac{1}{n}\left|\cro{\mu,\chi\Delta\phi} \sum_{\substack{k\in\N^*;\\k\le \cro{\mu,\chi}/4}}{\mu(k)} \sum_{\substack{y \in \M_F(\N);\\\cro{y,\mathbf 1}= k;\\\cro{y,\chi}\le n\alpha/2}}\P_Y(y \mid \mu,k, k) 
{ \cro{y,\chi}-k \over \cro{\mu,\chi}-\cro{y,\chi}-k} -{ \cro{y,\chi}- k \over \cro{\mu,\chi}} \right|} \\
+   \frac{1}{n}\left| \cro{\mu,\chi\Delta\phi}  \sum_{\substack{k\in\N^*;\\k\le \cro{\mu,\chi}/4}}{\mu(k)} \sum_{\substack{y \in \M_F(\N);\\\cro{y,\mathbf 1}=k;\\\cro{y,\chi}\le n\alpha/2}}\P_Y(y \mid \mu,k, k) 
{ \cro{y,\chi}-k \over \cro{\mu,\chi}} -\left({\cro{\mu,\chi^2} \over \cro{\mu,\chi}}-1\right)  \right|\\
+ \frac{2 \parallel \phi \parallel}{n}\left| \sum_{\substack{k\in\N^*;\\k\le \cro{\mu,\chi}/4}}{\mu(k)}\sum_{\substack{y \in \M_F(\N);\\\cro{y,\mathbf 1}= k;\\\cro{y,\chi}\le n\alpha/2}}\P_Y(y \mid \mu,k, k) 
\left(\cro{y,\chi}- k \right)\frac{ k+\cro{y,\chi}}{ \cro{\mu,\chi}-\cro{y,\chi}-k} 
\right|\\
:= \cro{\bar \B^n_{111}(\mu),\phi}+ \cro{\bar \B^n_{112}(\mu),\phi}+\cro{\bar \B^n_{113}(\mu),\phi}.
\end{multline}

We first have that
\begin{multline}
\begin{aligned}
\cro{\bar \B^n_{111}(\mu),\phi}&\leq 
2 M \parallel \phi \parallel\left|\sum_{\substack{k\in\N^*;\\k\le \cro{\mu,\chi}/4}}{\mu(k)}\sum_{\substack{y \in \M_F(\N);\\\cro{y,\mathbf 1}= k;\\\cro{y,\chi}\le n\alpha/2}}\P_Y(y \mid \mu,k,k) 
{ \cro{y,\chi}^2 \over \cro{\mu,\chi}(\cro{\mu,\chi}-\cro{y,\chi}-k)} \right|\\ 
&\quad\quad\quad+ 2M \parallel \phi \parallel\left|\sum_{\substack{k\in\N^*;\\k\le \cro{\mu,\chi}/4}}{\mu(k)}\sum_{\substack{y \in \M_F(\N);\\\cro{y,\mathbf 1}= k;\\\cro{y,\chi}\le n\alpha/2}}\P_Y(y \mid \mu,k,k) 
{ k^2 \over \cro{\mu,\chi}(\cro{\mu,\chi}-\cro{y,\chi}-k)} \right|\\
 &=o(n),\label{eq:convBn111}
 \end{aligned}
\end{multline}
since
\begin{equation*}
\left|\sum_{\substack{k\in\N^*;\\k\le \cro{\mu,\chi}/4}}{\mu(k)}\sum_{\substack{y \in \M_F(\N);\,\cro{y,\mathbf 1}= k;\\\cro{y,\chi}\le n\alpha/2}}\P_Y(y \mid \mu,k,k) 
{ k^2 \over \cro{\mu,\chi}(\cro{\mu,\chi}-\cro{y,\chi}-k)} \right|
 \leq  \frac{4}{(n\alpha)^2}\cro{\mu,\chi^2} \leq \frac{4}{n\alpha^2} M
\end{equation*}
and 
\begin{multline*}
\left|\sum_{\substack{k\in\N^*;\\k\le \cro{\mu,\chi}/4}}{\mu(k)}\sum_{\substack{y \in \M_F(\N);\,\cro{y,\mathbf 1}= k;\\\cro{y,\chi}\le n\alpha/2}}\P_Y(y \mid \mu,k,k) 
{ \cro{y,\chi}^2 \over \cro{\mu,\chi}(\cro{\mu,\chi}-\cro{y,\chi}-k)} \right|\ \\ 
\leq
\frac{4}{(n\alpha)^2}\sum_{\substack{k\in\N^*;\\k\le \cro{\mu,\chi}/4}}{\mu(k)} \E(\cro{Y(\mu),\chi}^2)
\leq \frac{4}{(n\alpha)^2} \frac{16 M^2}{3\alpha(3\alpha-4/n)} \sum_{\substack{k\in\N^*;\\k\le \cro{\mu,\chi}/4}}{\mu(k)} k^2 
\leq \frac{4}{n\alpha^2} \frac{16 M^3}{3\alpha(3\alpha-4/n)}, 
\end{multline*}
where we also used \eqref{eq:esperanceL2}. 
We now focus on the second term of \eqref{eq:convBn11}. From (\ref{eq:esperanceLbis}) we have that 
\begin{multline}
\label{eq:convBn112}
\cro{\bar \B^n_{112}(\mu),\phi}\\
\shoveleft{\le 2M \parallel \phi \parallel \left| \sum_{\substack{k\in\N^*;\\k\le \cro{\mu,\chi}/4}}{\mu(k)\over{\cro{\mu,\chi}}} 
 \left(k\left({\cro{\mu,\chi^2} \over \cro{\mu,\chi}-k}-1\right)+k^3o_3(n)+k(k-1)o_4(n)\right) -\left({\cro{\mu,\chi^2} \over \cro{\mu,\chi}}-1\right)  \right|}\\
\shoveleft{\leq 2 M \parallel \phi \parallel\left( \sum_{\substack{k\in\N^*;\\k\le \cro{\mu,\chi}/4}}{k\mu(k)\over{\cro{\mu,\chi}}} \left|
{\cro{\mu,\chi^2} \over \cro{\mu,\chi}-k}-{\cro{\mu,\chi^2} \over \cro{\mu,\chi}}  \right| +  \frac{1}{ \alpha}M  o_3(n) +  \frac{1}{ \alpha}M o_4(n)\right.} \\ \shoveright{\left.
+ \sum_{\substack{k\in\N^*;\\k> \cro{\mu,\chi}/4}}{k\mu(k)\over{\cro{\mu,\chi}}} \left|{\cro{\mu,\chi^2} \over \cro{\mu,\chi}}-1\right|\right)}\\
\leq 2 M \parallel \phi \parallel\left( \sum_{\substack{k\in\N^*;\\k\le \cro{\mu,\chi}/4}}{k\mu(k)\over{\cro{\mu,\chi}}} {k\cro{\mu,\chi^2}\over \cro{\mu,\chi}(\cro{\mu,\chi}-k)}+ \frac{1}{ \alpha}M (o_3(n)+o_4(n)) + \frac{M}{(n\alpha^2)}\frac{4}{(n\alpha)^2}M\right) \\
\leq   2 M \parallel \phi \parallel\left(  \frac{4M^2}{3n\alpha^3} +  \frac{1}{ \alpha}M(o_3(n)+o_4(n)) + \frac{4M^2}{n\alpha^3}\right),
\end{multline}
where we have used  the same argument as for (\ref{eq:convCn1}) for the last term. 
Finally,
\begin{multline*}
 \cro{\bar \B^n_{113}(\mu),\phi}
\leq \frac{2\parallel \phi \parallel}{n}  \left| \sum_{\substack{k\in\N^*;\\k\le \cro{\mu,\chi}/4}} \mu(k) \sum_{\substack{y \in \M_F(\N);\\\cro{y,\mathbf 1}= k;\\\cro{y,\chi}\le n\alpha/2}}\P_Y(y \mid \mu,k, k)  \frac{\cro{y,\chi}^2 -k^2}{\cro{\mu,\chi}-\cro{y,\chi}-k} \right|\\
\leq \frac{2\parallel \phi \parallel}{n}  \left( \frac{4}{n\alpha} \frac{16M^2}{3\alpha(3\alpha-4/n)} + \frac{4M}{\alpha}\right),
\end{multline*}
which also follows from \eqref{eq:esperanceL2}. This, together 
with (\ref{eq:convBn111}) and (\ref{eq:convBn112}) in (\ref{eq:convBn11}), concludes the proof.
\ep